\documentclass[12pt,a4paper]{amsart}

\usepackage[margin=1in]{geometry}  
\usepackage{graphicx}              
\usepackage{amsmath,euscript,enumerate}
\usepackage{amsfonts, color}              
\usepackage{amssymb, amscd,epsf,amsthm, epsfig, }                
\usepackage[toc,page]{appendix}

\newtheorem{thm}{Theorem}[section]

\newtheorem{prop}[thm]{Proposition}
\newtheorem{cor}[thm]{Corollary}
\newtheorem{lemma}[thm]{Lemma}
\newtheorem{claim}[thm]{Claim}
\theoremstyle{definition}
\newtheorem{rem}[thm]{Remark}

\input xy
\xyoption{all}

\usepackage{lipsum}

\newcommand{\Z}{\mathbb{Z}}
\newcommand{\LL}{\mathbb{L}}

\newcommand{\calF}{\mathcal{F}}

\newcommand{\bbL}{\mathbb{L}}

\newcommand{\bfn}{\mathbf{n}}

\newcommand{\DD}{\EuScript D}
\newcommand{\PP}{\EuScript P}

\newcommand{\KK}{\EuScript K}

\newcommand{\Hom}{\mathop{\mathrm{Hom}}\nolimits}

\newcommand{\Imm}{\mathop{\mathrm{Im}}\nolimits}

\newcommand{\sgDefo}{\mathop{\mathrm{sgDefo}}\nolimits}
\newcommand{\sgDPic}{\mathop{\mathrm{sgDPic}}\nolimits}

\newcommand{\Aut}{\mathop{\mathrm{Aut}}\nolimits}
\newcommand{\Barr}{\mathop{\mathrm{Bar}}\nolimits}
\newcommand{\modu}{\mathop{\mathrm{mod}}\nolimits}
\newcommand{\Modu}{\mathop{\mathrm{Mod}}\nolimits}

\newcommand{\perf}{\mathop{\mathrm{Perf}}\nolimits}
\newcommand{\id}{\mathop{\mathrm{Id}}\nolimits}

\newcommand{\op}{\mathop{\mathrm{op}}\nolimits}

\newcommand{\Ker}{\mathop{\mathrm{Ker}}\nolimits}

\newcommand{\sg}{\mathop{\mathrm{sg}}\nolimits}
\newcommand{\sy}{\mathop{\mathrm{sy}}\nolimits}
\newcommand{\nc}{\mathop{\mathrm{nc}}\nolimits}

\newcommand{\HH}{\mathop{\mathrm{HH}}\nolimits}

\newcommand{\Lie}{\mathop{\mathrm{Lie}}\nolimits}

\newcommand{\RHS}{\mathop{\mathrm{RHS}}\nolimits}

\makeatletter
\newcommand{\colim@}[2]{%
  \vtop{\m@th\ialign{##\cr
    \hfil$#1\operator@font colim$\hfil\cr
    \noalign{\nointerlineskip\kern1.5\ex@}#2\cr
    \noalign{\nointerlineskip\kern-\ex@}\cr}}%
}
\newcommand{\colim}{%
  \mathop{\mathpalette\colim@{\rightarrowfill@\textstyle}}\nmlimits@
}
\makeatother

\begin{document}


\title[Invariance of the Gerstenhaber structure on Tate-Hochschild]{Invariance of the Gerstenhaber algebra structure on Tate-Hochschild cohomology}
\date{}

\author{Zhengfang WANG}
 \address{Beijing International Center for Mathematical Research (BICMR), Peking University, No. 5 Yiheyuan Road Haidian District, Beijing 100871, P.R. China} \address{Universit\'e Paris Diderot-Paris 7, Institut de Math\'ematiques de Jussieu-Paris Rive Gauche CNRS UMR 7586, B\^atiment Sophie Germain, Case 7012, 75205 Paris Cedex 13, France}

\email{wangzhengfang@bicmr.pku.edu.cn}


\maketitle
\begin{abstract}
Keller proved in 1999 that the Gerstenhaber algebra structure on the Hochschild cohomology of an algebra is an invariant of the derived category. In this paper, we adapt his approach to show that the  Gerstenhaber algebra structure on the Tate-Hochschild cohomology of an algebra is preserved under singular equivalences of Morita type with level, a notion  introduced by the author in  previous work.

\end{abstract}

\section{Introduction}
In \cite{Wan15a, Wan15}, we constructed a Gerstenhaber algebra structure on the Tate-Hochschild cohomology ring $\HH^*_{\sg}(A, A)$ implicit in Buchweitz' work \cite{Buc} for an algebra $A$ projective over a commutative ring $k$ and such that $A$ and the enveloping algebra $A\otimes_k A^{\op}$ are Noetherian. The cup product is given by the Yoneda product in the singularity category of the enveloping algebra $A\otimes_k A^{\op}$.  Recall that the {\it singularity category} $\DD_{\sg}(A)$ (cf. \cite{Buc, Orl}) of a Noetherian algebra $A$ is defined  as the Verdier quotient of the bounded derived category $\DD^b(A)$ of finitely generated (left) $A$-modules by the full subcategory $\perf(A)$ consisting of complexes quasi-isomorphic to bounded complexes of finitely generated projective $A$-modules. The Lie bracket  on $\HH_{\sg}^*(A, A)$ was defined in \cite{Wan15a, Wan15} as the graded commutator of a certain  circle product $\circ$ extending naturally the Gerstenhaber circle product on Hochschild cohomology.   In particular, for a self-injective algebra, in positive degrees, this Lie bracket coincides with the Gerstenhaber bracket in Hochschild cohomology.  In \cite{Wan15a, Wan15}, we also proved that the natural morphism, induced by the quotient functor from the bounded derived category to the singularity category of $A\otimes_k A^{\op}$,   from the Hochschild cohomology ring $\HH^*(A, A)$ to $\HH^*_{\sg}(A, A)$ is a morphism of Gerstenhaber algebras. By the very recent work of Keller \cite{Kel18}, the Tate-Hochschild cohomology of an algebra $A$ is isomorphic, as graded algebras,  to the Hochschild cohomology of the dg singularity category (i.e. the canonical dg enhancement of the singularity category) of $A$. This yields a second Gerstenhaber algebra structure on Tate-Hochschild cohomology, which is conjectured    to coincide with the one introduced in \cite{Wan15a, Wan15}. For more details, we refer to Keller's conjecture \cite[Conjecture 1.2]{Kel18}.

Keller proved in \cite{Kel} that the Gerstenhaber algebra structures on Hochschild cohomology rings  are preserved under derived equivalences of standard type. That is, let   $X$ be a complex of $A$-$B$-bimodules such that the total derived tensor product by $X$ is an equivalence  between the derived categories of two $k$-algebras $A$ and $B$. Then $X$ yields a natural isomorphism of Gerstenhaber algebras from $\HH^*(A, A)$ to $\HH^*(B,B)$. In this paper, we will show that the Gerstenhaber algebra structure on the Tate-Hochschild cohomology ring is also preserved under derived equivalences of standard type. In fact,  we will prove a stronger result.  Namely,   the Gerstenhaber algebra structure on the Tate-Hochschild cohomology ring is preserved under singular equivalences of Morita type with level (cf. \cite{Wan14} and Section \ref{section6} below). Recall that a derived equivalence of standard type induces a singular equivalence of Morita type with level (cf. \cite{Wan14}).

The paper  is organized as follows. In Section 2, we recall the construction  of the normalized bar resolution of an algebra $A$ and provide some natural liftings of elements in $\HH^*_{\sg}(A, A)$ along the normalized bar resolution.  In Section 3,  we introduce the bullet product $\bullet$ and the circle product $\circ$. Using these two products,  we construct two dg modules $C^L(f, g)$ and $C^R(f, g)$ associated to the cohomology classes $f$ and $g$ in $\HH^*(A, \Omega_{\sy}^*(A))$. These two dg modules  play a crucial role in the proof of  our main result.  In Section 4, we  recall the notions of $R$-relative derived categories and $R$-relative derived tensor products. In Section 5, we  develop the singular infinitesimal deformation theory of the identity bimodule in analogy with the infinitesimal deformation theory of  \cite{Kel}. As a result, we give an interpretation of the Gerstenhaber bracket on the Tate-Hochschild cohomology ring from the point of view of the singular infinitesimal deformation theory. In Section 6, we prove our main result.
\begin{thm}[=Theorem \ref{thm} and Corollary \ref{cor-main}]
Let $k$ be a field. Let $A$ and $B$ be two Noetherian (not necessarily commutative) $k$-algebras such that the enveloping algebras $A\otimes A^{\op}$ and $B\otimes B^{\op}$ are Noetherian. Suppose that $(_AM_B, _BN_A)$ defines a singular equivalence of Morita type with level $l\in\Z_{\geq 0}$. Then the functor $$\Sigma^l (M\otimes_B-\otimes_BN): \DD_{\sg}(B\otimes_k B^{\op})\rightarrow \DD_{\sg}(A\otimes_k A^{\op})$$ induces an isomorphism of Gerstenhaber algebras between the Tate-Hochschild cohomology rings $\HH_{\sg}^*(A, A)$ and $\HH_{\sg}^*(B, B)$. In particular, the Gerstenhaber algebra structure on the Tate-Hochschild cohomology ring is invariant under derived equivalences. 
\end{thm}

\begin{rem}
Let $k$ be an algebraically closed field. Let $A$ and $B$ be two (finite dimensional) symmetric $k$-algebras which are related by a stable equivalence of Morita type. Then the authors in \cite[Theorem 10. 7]{KLZ} proved that there is an isomorphism of Gerstenhaber algebras (more generally, BV algebras)  between $\HH_{\sg}^{\geq 0}(A, A)$ and $\HH_{\sg}^{\geq 0}(B, B)$. 
\end{rem}

Throughout this paper, we fix a field $k$. The unadorned tensor product  $\otimes$ and $\Hom$ represent the tensor product $\otimes_k$ and  $\Hom_k$ over the field $k$, respectively. We write the composition $g\circ f$ of two maps $f: X\rightarrow Y$ and $g: Y\rightarrow Z$ as $gf$. We write  the identity map  $\id_X: X\rightarrow X$ simply as $\id$ when no confusion can arise.  We will follow the Koszul sign rule for the tensor product: $(f\otimes g)(x\otimes y)=(-1)^{|g||x|} f(x)\otimes g(y)$ where $|g|$ is the degree of the homogeneous map $g$ and $|x|$ is the degree of the element $x\in X$.    

The notions of differential graded (dg) algebras and relative tensor products are frequently used in this paper.  For more details, we refer to \cite{Kel, Kel2, BeLu},  and to \cite{KeVo, Ric1, Wei, Zim} for the notions of triangulated categories and derived categories.

\section*{Acknowledgement} This paper is part of the author's Ph. D.  thesis, the author would like to thank his  supervisor Alexander Zimmermann for introducing him to  this interesting topic and for many valuable suggestions for improvement. He is very grateful to Bernhard Keller for many useful discussions and suggestions on relative derived categories and relative singularity categories.  He  would like to thank Gufang Zhao for many useful  discussions at the very beginning  of this project. He is also very grateful to Yu Wang for reading a preliminary version of  this paper and for giving many helpful  suggestions. 

The author would like to thank the referee for valuable suggestions and comments, which have led to significant improvement on the presentation of this paper.

\section{Normalized bar resolution}
\subsection{The definition}\label{section2.1}
Let $A$ be an associative algebra over a field $k$. The  {\it normalized bar resolution}  (cf. e.g. \cite{Lod}) is defined as the dg  $A$-$A$-bimodule $\Barr_*(A):=\bigoplus_{p\geq 0} \Barr_p(A),$ with $\Barr_p(A):=A\otimes (\Sigma \overline{A})^{\otimes p}\otimes A \ (p\geq 0)$ in degree $p$ and the differential of degree $-1$ 
 \begin{equation*}
\begin{split}
d_p(a_0\otimes \overline{a}_{1, p}\otimes a_{p+1})=&a_0a_1\otimes \overline{a}_{2, p}\otimes a_{p+1}+\\
&\sum_{i=1}^{p-1} (-1)^ia_0\otimes \overline{a}_{1, i-1}\otimes \overline{a_ia_{i+1}}\otimes \overline{a}_{i+2, p}\otimes a_{p+1}+\\
& (-1)^{p}a_0\otimes \overline{a}_{1, p-1}\otimes a_pa_{p+1}.
\end{split}
\end{equation*}
Let us explain the notations appeared above:  We denote by $\Sigma \overline{A}$ the graded $k$-module concentrated in degree $1$ with $(\Sigma \overline{A})_1=A/(k\cdot 1); $ Let $\pi: A\rightarrow (\Sigma \overline{A})$ be the natural projection of degree $1$. Then we denote $\overline a=\pi(a)$  for any $a\in A$. The degree of $\overline a$ is $|\overline{a}|=1$;  We simply  write   $\overline{a_i}\otimes \overline{a_{i+1}}\otimes \cdots \otimes \overline{a_j}\in(\Sigma \overline{A})^{\otimes (j-i+1)}$ as  $\overline{a}_{i, j}$.   It is well-known that $\Barr_*(A)$ is a projective  bimodule resolution of  $A$ with the augmentation map $\tau_0=d_0: A\otimes A \rightarrow A, a\otimes b \mapsto ab$. For convenience, we set $\Barr_{-1}(A)=A$.  

For any $p\in \Z_{>0}$, we denote the kernel of
the differential $d_{p-1}: \Barr_{p-1}(A) \rightarrow \Barr_{p-2}(A)$ by $\Omega_{\sy}^p(A)$. In particular, we set $\Omega_{\sy}^0(A)=A$. It is clear that $\Omega_{\sy}^p(A)$ is an $A$-$A$-bimodule. For convenience, we view $\Omega_{\sy}^p(A)$ as a dg bimodule concentrated in degree $p$.  
For $p\geq 0$, we denote by   $\Barr_{\geq p}(A)$ the  \lq $p$-truncated'  dg $A$-$A$-bimodule  with $\Barr_{\geq p}(A)_i=\Barr_i(A)$ if $i\geq p$ and $\Barr_{\geq p}(A)_i=0$ if $i<p$.  Recall that, for a chain complex $(X, d)$ and  $p\in\Z$,   the {\it $p$-shifted complex}  $(\Sigma^pX, \Sigma^pd)$ is defined as $(\Sigma^pX)_n=X_{n-p}$  with the differential $(\Sigma^pd)_n=(-1)^p d_{n-p}$ for any $n\in \mathbb Z.$

\begin{rem}\label{rem-homotopy}
Note that the  \lq $p$-truncated' augmented  normalized bar resolution     $$\widehat{\Barr}_{\geq p}(A):\cdots \rightarrow \Barr_{p+1}(A)\xrightarrow{d_{p+1}} \Barr_p(A)\xrightarrow{\tau_p=d_p}  \Sigma^{-1}\Omega_{\sy}^p(A)\rightarrow 0 $$ is exact for any fixed $p\in\Z_{\geq 0}$.  For this, we define a $k$-linear map for any $r\geq 0$,
\begin{equation*}
s^L_r: \Barr_r(A)  \rightarrow \Barr_{r+1}(A), \quad a_0\otimes \overline{a}_{1, r}\otimes a_{r+1} \mapsto (-1)^{r+1} a_0\otimes  \overline{a}_{1, r+1}\otimes 1. \end{equation*}
\begin{equation*}
\xymatrix@C=4pc@R=3.4pc{
\cdots\ar[r]  & \Barr_{p+1}(A)\ar[d]^-{\id} \ar[r]^-{d_{p+1}}  & \Barr_{p}(A)\ar[ld]^-{s^L_{p}} \ar[r]^-{d_p} \ar[d]_-{\id} &   \Sigma^{-1}\Omega_{\sy}^p(A)\ar[r] \ar[d]^-{\id}  \ar[ld]^-{s^L_{p-1}|_{\Sigma^{-1}\Omega_{\sy}^p(A)}} & 0\\
\cdots\ar[r]  &  \Barr_{p+1}(A)  \ar[r]_-{d_{p+1}}&  \Barr_{p}(A) \ar[r]_-{d_p}  &  \Sigma^{-1} \Omega_{\sy}^p(A) \ar[r]   & 0
}
\end{equation*}
It is straightforward to verify that $s^L  d+d  s^L=\id_{\widehat{\Barr}_{\geq p}(A)}.$ This yields the exactness of $\widehat{\Barr}_{\geq p}(A)$. Note that $s^L$ is a morphism of  left graded $A$-modules (but  not a morphism of graded $A\otimes A^{\op}$-modules). 
Similarly, if we define 
$$s_r^R: \Barr_r(A)\rightarrow \Barr_{r+1}(A),\quad a_0\otimes \overline{a}_{1, r}\otimes a_{r+1}\mapsto 1\otimes \overline{a}_{0, r}\otimes a_{r+1}, $$
then we have that $s^R d+d s^R=\id_{\widehat{\Barr}_{\geq p}(A)}$ and $s^R$ is a morphism of right graded  $A$-modules.  
\end{rem}

For any $p, q\in \Z_{\geq 0}$, we will construct a morphism of
dg $A$-$A$-bimodules between 
$\Barr_{\geq p+q}(A)$ and $\Barr_{\geq p}(A)\otimes_A\Barr_{\geq q}(A)$.
We define
\begin{equation*}
  \Delta_{p, q}: \Barr_{\geq p+q}(A)\rightarrow \Barr_{\geq p}(A)
  \otimes_A\Barr_{\geq q}(A)
\end{equation*}
as follows.  For $a_0\otimes \overline{a}_{1, p+q+r}\otimes a_{p+q+r+1}\in
\Barr_{p+q+r}(A)$, where $ r\geq 0$, 
$$\Delta_{p, q }(a_0\otimes \overline{a}_{1, p+q+r}\otimes a_{p+q+r+1})=\sum_{i=0}^r
(a_0\otimes \overline{a}_{1, p+i}\otimes 1)\otimes_A( 1\otimes \overline{a}_{p+i+1, p+q+r}\otimes a_{p+q+r+1}).$$
 It is a routine computation to verify that $\Delta_{p, q}$ is a morphism of dg $A$-$A$-bimodules.

\begin{lemma}\label{lemma2.2}
For any  $p,q\in\Z_{\geq 0}$,  $\Delta_{p, q}$ is an isomorphism in the homotopy category $\KK(\mbox{$A\otimes A^{\op}$-$\Modu$})$ of dg $A$-$A$-bimodules.
\end{lemma}
\begin{proof} For $p=0$ or $q=0$, we have a natural  isomorphism $\mu_{p, q}: \Omega_{\sy}^p(A)\otimes_A \Omega_{\sy}^q(A)\xrightarrow{\cong} \Omega_{\sy}^{p+q}(A)$ since $\Omega_{\sy}^0(A)=A$. 
For $p, q>0$, consider  the following composition of maps 
$$\mu_{p, q}: \Omega_{\sy}^p(A)\otimes_A\Omega_{\sy}^q(A)\hookrightarrow  A\otimes (\Sigma \overline{A})^{\otimes p-1}\otimes  A\otimes (\Sigma \overline{A})^{\otimes q-1}\otimes A\xrightarrow{\id^{\otimes p}\otimes \pi\otimes \id^{\otimes q}} A\otimes (\Sigma \overline{A})^{\otimes p+q-1} \otimes A,$$ where the first map is given by the tensor product of the natural inclusions 
$$\Omega_{\sy}^p(A)\hookrightarrow \Barr_{p-1}(A),\quad \Omega_{\sy}^q(A)\hookrightarrow \Barr_{q-1}(A),$$ and where  $\pi: A\rightarrow \Sigma \overline{A}$ is the natural projection of degree $1$.  More concretely, let $x:=\sum_i a_0^i\otimes \overline{a^i}_{1, p-1}\otimes a_{p}^i\in \Omega_{\sy}^{p}(A)$ and $y:=\sum_j b_0^j\otimes \overline{b^j}_{1, q-1}\otimes b_{q}^j\in \Omega_{\sy}^{q}(A)$.  Then 
$$\mu_{p, q}(x\otimes_Ay)=\sum_{i, j} a_0^i\otimes \overline{a^i}_{1, p-1}\otimes \overline{a_{p}^i b_0^j}\otimes \overline{b^j}_{1, q-1}\otimes b_{q}^j.$$ Notice that the image of $\mu_{p, q}$ lies in $\Omega_{\sy}^{p+q}(A)$ since $d_{p+q-1}\mu_{p, q}(x\otimes_A y)=0$.  This induces an $A$-$A$-bimodule homomorphism $\mu_{p, q}:  \Omega_{\sy}^p(A)\otimes_A\Omega_{\sy}^q(A)\rightarrow \Omega_{\sy}^{p+q}(A)$.     We claim  that  $\mu_{p, q}$ is a bijection and its inverse  $\mu_{p, q}^{-1}: \Omega^{p+q}_{\sy}(A)\rightarrow \Omega_{\sy}^p(A)\otimes_A \Omega_{\sy}^q(A)$  sends an element $x:=\sum_i a_0^i\otimes \overline{a^i}_{1, p+q-1}\otimes a_{p+q}^i\in \Omega_{\sy}^{p+q}(A)$ to   
$$\mu_{p, q}^{-1}(x)=(-1)^{p+q}\sum_i d_p(a_0^i\otimes \overline{a^i}_{1, p}\otimes 1)\otimes_A d_q(1\otimes \overline{a^i}_{p+1, p+q}\otimes 1).$$ Indeed,  we  have 
\begin{equation*}
\begin{split}
\mu_{p, q}\mu^{-1}_{p, q}(x)=&(-1)^{p+q} \sum_i a_0^ia_1^i \otimes \overline{a^i}_{2, p+q}\otimes 1+\\
&\sum_i\sum^{p+q-1}_{j=1} (-1)^{p+q+j} a_0^i\otimes \overline{a^i}_{1, j-1} \otimes \overline{a^i_ja^i_{j+1}}\otimes \overline{a^i}_{j+2, p+q}\otimes 1+\\
& \sum_i a_0^i\otimes \overline{a^i}_{1, p+q-1}\otimes a_{p+q}^i\\
=&(\id^{\otimes p+q-2}\otimes \pi)(dx)\otimes 1+x\\
=&x,
 \end{split}
\end{equation*} 
where the third identity comes from the identity $dx=0$ (since $x\in \Omega_{\sy}^{p+q}(A)$).  Similarly, for $x:=\sum_i a_0^i \otimes \overline{a^i}_{1, p-1}\otimes a^i_p\in \Omega^p_{\sy}(A)$ and $y:=\sum_j b_0^j\otimes \overline{b^j}_{1, q-1}\otimes b^j_q\in \Omega^q_{\sy}(A)$, we have 
\begin{equation*}
\begin{split}
\mu^{-1}_{p, q}\mu_{p, q}(x\otimes_A y)=&\mu^{-1}_{p, q} \left(\sum_{i, j} a_0^i \otimes \overline{a^i}_{1, p-1}\otimes \overline{a^i_pb^j_0}\otimes \overline{b^j}_{1, q-1}\otimes b_q^j\right)\\
=&(-1)^{p+q} \sum_{i, j} d_p(a_0^i \otimes \overline{a^i}_{1, p-1}\otimes \overline{a^i_pb_0^j}\otimes 1)\otimes_A d_q(1\otimes \overline{b^j}_{1, q}\otimes 1)\\
=& x\otimes_A y\end{split}
\end{equation*}
where the third identity comes from the identities $dx=0$ and $dy=0$. This proves the claim. It is clear that   $\mu_{p, q}$ is a morphism of $A$-$A$-bimodules. Hence  so is  $\mu_{p, q}^{-1}$. 
Since $(\tau_p\otimes_A \tau_q) \Delta_{p, q}= \mu_{p, q}^{-1}\tau_{p+q}$, we get that  $\Delta_{p, q}$ is a lifting of the isomorphism $\mu^{-1}_{p, q}$ between the  resolutions $\Barr_{\geq p+q}(A)$ and $\Barr_{\geq p}(A)\otimes_A\Barr_{\geq q}(A)$.  Hence  it  is an isomorphism in the homotopy category $\KK(\mbox{$A\otimes A^{\op}$-$\Modu$})$ of dg $A$-$A$-bimodules. \end{proof}


For $p\geq 0$, we define  the dg $A$-$A$-bimodule of left  noncommutative differential $p$-forms as $\Omega_{\nc}^{L, p}(A) =A\otimes (\Sigma \overline{A})^{\otimes p}$. Clearly, $\Omega_{\nc}^{L, p}(A)$ is concentrated in degree $p$.  The  bimodule structure  is given by 
$$a(a_0\otimes \overline{a}_{1, p})\blacktriangleleft b=-(\id^{\otimes p}\otimes \pi) d_p(aa_0\otimes \overline{a}_{1, p}\otimes b) $$
for $a, b \in A$ and $a_0\otimes \overline{a}_{1, p}\in A\otimes (\Sigma \overline{A})^{\otimes p}$. Here when  $(\id^{\otimes p}\otimes \pi)$ is applied  to the element $d_p(aa_0\otimes \overline{a}_{1, p}\otimes b)$, additional signs will appear because of the Koszul sign rule since $\pi$ is a map of degree $1$. More explicitly, we have
\begin{equation*}
\begin{split}
a(a_0\otimes \overline{a}_{1, p})\blacktriangleleft b=&(-1)^{p} aa_0a_1\otimes \overline{a}_{2, p}\otimes \overline{b}+\\
&\sum_{i=1}^{p-1} (-1)^{p+i} aa_0\otimes \overline{a}_{1, i-1}\otimes \overline{a_ia_{i+1}}\otimes \overline{a}_{i+2, p}\otimes \overline{b}\\&+aa_0\otimes\overline{a}_{1, p-1}\otimes \overline{a_pb}.
\end{split}
\end{equation*} Similarly, the dg $A$-$A$-bimodule of right  noncommutative differential $p$-forms is defined  as $\Omega_{\nc}^{R, p}(A)=(\Sigma \overline{A})^{\otimes p}\otimes A$. The bimodule structure is given by 
$$a\blacktriangleright (\overline{a}_{1, p}\otimes a_{p+1})b=(\pi\otimes \id^{\otimes p}) d_p(a\otimes \overline a_{1, p}\otimes a_{p+1}b). $$

The following lemma is very useful throughout the present paper. 
\begin{lemma}\label{lemma-forms}
We have two isomorphisms of dg $A$-$A$-bimodules
$$\alpha_p^L: \Omega_{\nc}^{L, p}(A)\xrightarrow{\cong} \Omega_{\sy}^p(A),\quad a_0\otimes \overline{a}_{1, p}\mapsto -d_p(a_0\otimes \overline{a}_{1, p}\otimes 1) ;$$
$$ \alpha_p^R:\Omega_{\nc}^{R, p}(A) \xrightarrow{\cong} \Omega_{\sy}^p(A), \quad \overline{a}_{1, p}\otimes a_0\mapsto d_p(1\otimes \overline{a}_{1, p}\otimes a_0).$$
\end{lemma}
\begin{proof}
First, we claim that both $\alpha_p^L$ and $\alpha_p^R$ are bijective. Indeed, the inverse of $\alpha_p^L$ is given by 
$$(\alpha_p^L)^{-1}(x)=(-1)^{p-1}\sum_{i} a_0^i\otimes \overline{a^i}_{1, p}$$
for $x:=\sum_i a_0^i \otimes \overline{a^i}_{1, p-1}\otimes a^i_p\in \Omega^p_{\sy}(A)$.  That is, $(\alpha_p^L)^{-1}$ is the composition of maps
$$(\alpha_p^L)^{-1}: \Omega^p_{\sy}(A)\hookrightarrow A\otimes (\Sigma \overline{A})^{\otimes p-1}\otimes A\xrightarrow{\id^{\otimes p}\otimes \pi} \Omega^{L, p}_{\nc}(A).$$
Here the sign $(-1)^{p-1}$ is hidden in  the Koszul sign rule. 
From a straightforward computation, we get that 
$\alpha_p^L(\alpha_p^L)^{-1}=\id$ and $(\alpha_p^L)^{-1}\alpha_p^L=\id.$ Similarly, the inverse of $\alpha_p^R$ is given by 
$$(\alpha_p^R)^{-1}(x)=\sum_{i}  \overline{a^i}_{0, p-1}\otimes a^i_p$$
for $x:=\sum_i a_0^i \otimes \overline{a^i}_{1, p-1}\otimes a^i_p\in \Omega^p_{\sy}(A)$. That is, $(\alpha_p^R)^{-1}$ is the composition of maps
$$(\alpha_p^R)^{-1}: \Omega_{\sy}^p(A)\hookrightarrow A\otimes (\Sigma \overline{A})^{\otimes p-1}\otimes A\xrightarrow{\pi\otimes \id^{\otimes p}} \Omega_{\nc}^{R, p}(A).$$ This proves the claim. It remains to check that $\alpha_p^L$ and $\alpha_p^R$ are morphisms of $A$-$A$-bimodules. For this, given  $a_0\otimes \overline{a}_{1, p}\in A\otimes (\Sigma \overline{A})^{\otimes p}$, we have 
\begin{equation*}
\begin{split}
\alpha_p^L(a(a_0\otimes\overline{a}_{1, p})\blacktriangleleft b)=&- d_p((\id^{\otimes p}\otimes \pi)d_p(aa_0\otimes \overline{a}_{1, p}\otimes b)\otimes 1)\\
=&d_p(aa_0\otimes \overline{a}_{1, p}\otimes b)\\
=& a\alpha_p^L(a_0\otimes \overline{a}_{1, p})b,
\end{split}
\end{equation*} 
where the second identity follows from $d_pd_{p+1}(aa_0\otimes\overline{a}_{1, p}\otimes \overline{b}\otimes 1)=0$. 
This shows that $\alpha_p^L$ is a morphism of $A$-$A$-bimodules. By a similar computation, we get that $\alpha_p^R$  is a morphism of $A$-$A$-bimodules. This proves the lemma. 
\end{proof}

\subsection{Two liftings }\label{section2.2}
Let $M$ be a graded  $A$-$A$-bimodule. Recall that the Hochschild cohomology $\HH^*(A, M)$ with coefficients in $M$ is computed by the Hochschild cochain complex $(C^*(A,M), \delta)$ with $$C^{m}(A, M)=\prod_{i\geq 0} \mathcal Hom^{-m}((\Sigma \overline{A})^{\otimes i}, M), \quad \mbox{for $m\in \Z$},$$ where  $(\Sigma\overline A)^{\otimes 0}=k$ and $\mathcal Hom^{-m}((\Sigma \overline{A})^{\otimes i}, M)$ is the set of $k$-linear maps of degree $-m$ from chain complexes $(\Sigma \overline{A})^{\otimes i}$ to $M$. Recall  that  a $k$-linear map $f: X\rightarrow Y$ between two chain complexes $X$ and $Y$ is {\it  of degree $m$} if  $f(X_i)\subset Y_{i+m}$ for any $i\in\mathbb Z$.   The differential $\delta$ (of degree one) is given by, for $f\in C^{m}(A, M)$,  
\begin{equation*}
\begin{split}
\delta^{m}(f)(\overline{a}_{1, i+1})=&a_1f(\overline{a}_{2, i+1})+\sum_{j=1}^i(-1)^j f(\overline{a}_{1, j-1}\otimes \overline{a_ja_{j+1}}\otimes \overline{a}_{j+2, i+1})+\\
& (-1)^{i+1}f(\overline{a}_{1, i})a_{i+1}.
\end{split}
\end{equation*}

Let $m, p\in \Z_{\geq 0}$ and $f\in \HH^{m-p}(A, \Omega_{\sy}^p(A))$. Recall that $\Omega_{sy}^p(A)$ is a graded $A$-$A$-bimodule concentrated in degree $p$.  Then $f$ can be represented by an element $f\in C^{m-p}(A, \Omega_{\sy}^p(A))=\Hom((\Sigma \overline{A})^{\otimes m}, \Omega_{\sy}^p(A))$ such that $\delta(f)=0.$ Denote 
$$f^L: (\Sigma \overline{A})^{\otimes m} \xrightarrow{f} \Omega_{\sy}^p(A)\xrightarrow{(\alpha^L_p)^{-1}} \Omega_{\nc}^{L, p}(A)=A\otimes (\Sigma \overline{A})^{\otimes p}; $$
$$f^R:(\Sigma \overline{A})^{\otimes m} \xrightarrow{f} \Omega_{\sy}^p(A)\xrightarrow{(\alpha^R_p)^{-1}} \Omega_{\nc}^{R, p}(A)= (\Sigma \overline{A})^{\otimes p}\otimes A,$$
where $(\alpha_p^L)^{-1}$ and $(\alpha_p^R)^{-1}$ are defined in Lemma \ref{lemma-forms}.  
  These two maps induce two liftings $$\vartheta^L(f), \vartheta^R(f): \Barr_*(A)\rightarrow \Sigma^{m-p}\Barr_{\geq p}(A)$$ in the following way.  Let  $x:=a_0\otimes \overline{a}_{1, r}\otimes a_{r+1}\in \Barr_{r}(A)$. If $r<m$,   we define $$\vartheta^L(f)(x)=\vartheta^R(f)(x)=0.$$ If $r\geq m$, we define 
\begin{equation*}
  \begin{split}
    \vartheta^L(f)(x)&=a_0f^L(\overline{a}_{1, m})\otimes \overline{a}_{m+1, r}\otimes a_{r+1},\\
   \vartheta^R(f)(x)&=(-1)^{(m-p)(r-m)}a_0\otimes \overline{a}_{1, r-m}\otimes f^R(\overline{a}_{r-m+1, r})a_{r+1}.
 \end{split}
\end{equation*}
 It follows from $\delta(f)=0$ that  $\vartheta^L(f)$ and $\vartheta^R(f)$ are indeed morphisms of dg $A$-$A$-bimodules. It is well-known from homological algebra (cf. e.g. \cite[Comparison Theorem 2.2.6]{Wei}) that $\vartheta^L(f)$ is homotopy equivalent to  $\vartheta^R(f)$. In fact, there exists a specific chain homotopy 
\begin{equation}\label{equ-chain-hom} 
h(f): \Barr_*(A) \rightarrow \Sigma^{m-p-1}\Barr_{\geq p}(A)
\end{equation} 
from $\vartheta^L(f)$ to $\vartheta^R(f)$ defined as follows. For any $r\in \Z_{\geq 0}$,  
\begin{eqnarray*}
\lefteqn{h_r(f)(a_0\otimes \overline{a}_{1, r}\otimes a_{r+1})} \\
&=\begin{cases}
0 & \mbox{for $r\leq m-1,$}\\
\sum\limits_{i=0}^{r-m} (-1)^{i(m-p-1)} a_0\otimes \overline{a}_{1, i} \otimes \overline{f}(\overline{a}_{i+1, i+m})\otimes \overline{a}_{i+m+1, r}\otimes a_{r+1} & \mbox{for $r\geq m,$}
\end{cases}
\end{eqnarray*}
where $$\overline{f}: (\Sigma \overline{A})^{\otimes m} \xrightarrow{f} \Omega_{\sy}^p(A)\hookrightarrow \Barr_{p-1}(A) \xrightarrow{\pi\otimes \id^{\otimes p-1}\otimes \pi} (\Sigma \overline{A})^{p+1}. $$
Indeed, it is easy to verify that $ \vartheta^L(f)-\vartheta^R(f)=h(f)d+d h(f).$ Notice that  $h(f)$ is a morphism of graded $A$-$A$-bimodules. It follows that $\vartheta^L(f)$ is isomorphic to $\vartheta^R(f)$ in the homotopy category $\KK^-(\mbox{$A\otimes A^{\op}$-$\Modu$})$. Therefore,  both $\vartheta^L(f)$ and $\vartheta^R(f)$ are representatives of   $f\in \Hom_{\DD^b(A\otimes A^{\op})}(A, \Sigma^{m-p}\Omega_{\sy}^p(A))$ in $\KK^-(\mbox{$A\otimes A^{\op}$-$\Modu$})$. 

From $f\in \HH^{m-p}(A, \Omega_{\sy}^p(A))$, we may get   an element $\Omega_{\sy}^r(f)\in \HH^{m-p}(A, \Omega_{\sy}^{p+r}(A)) $ for any $r\geq 0$, which is represented by the element 
\begin{equation*}
 \begin{tabular}{rccc}
 $ \Omega_{\sy}^{L, r}(f):$ & $\Barr_{m+r}(A)$ &$\rightarrow$& $ \Sigma^{m-p}\Omega_{\sy}^{p+r}(A)$\\
 & $a_0\otimes \overline{a}_{1, m+r}\otimes a_{m+r+1}$ & $\mapsto$ &  $d_{p+r}(a_0f^L(\overline{a}_{1, m})\otimes \overline{a}_{m+1, m+r}\otimes a_{m+r+1}).$
  \end{tabular}
\end{equation*}
Similarly, $\Omega_{\sy}^r(f)$ may  also be represented by the element
\begin{equation*}
 \begin{tabular}{rccc}
 $ \Omega_{\sy}^{R, r}(f):$ & $\Barr_{m+r}(A)$ &$\rightarrow$& $ \Sigma^{m-p}\Omega_{\sy}^{p+r}(A)$\\
 & $a_0\otimes \overline{a}_{1, m+r}\otimes a_{m+r+1}$ & $\mapsto$ &  $(-1)^{(m-p)r}d_{p+r}(a_0\otimes  \overline{a}_{1, r}\otimes f^R(\overline{a}_{r+1, m+r}) a_{m+r+1}).$
  \end{tabular}
\end{equation*}
\begin{rem}\label{remark-delta-cup}
The above homotopy $h(f)$ induces a homotopy $h^{L, R}_r(f):=d_{p+r}h_{m+r-1}(f)$ such that 
$h^{L, R}_r(f)d_{m+r}=\Omega_{\sy}^{L, r}(f)-\Omega_{\sy}^{R, r}(f).$
For any $f\in C^{m-p}(A, \Omega_{\sy}^p(A))$ such that $\delta(f)=0$, we have the following identities
$$\mu_{r, p+s}(d_r\otimes_A\Omega_{\sy}^{R, s}(f))\Delta_{r, m+s}=\Omega_{\sy}^{R, r+s}(f),$$
$$\mu_{p+s, r}(\Omega_{\sy}^{L, s}(f)\otimes_A d_r)\Delta_{m+s, r}=\Omega_{\sy}^{L, r+s}(f),$$
which can be verified by straightforward computation. 
\end{rem}
Therefore, we have a map for any $r>0$,  
$$\Omega^r_{\sy}:\HH^{m-p}(A, \Omega_{\sy}^p(A))\rightarrow \HH^{m-p}(A, \Omega_{\sy}^{p+r}(A)), \quad f\mapsto \Omega_{\sy}^{L, r}(f)=\Omega_{\sy}^{R, r}(f). $$
Notice that $\Omega^r_{\sy}(\Omega^s_{\sy}(f))=\Omega^{r+s}_{\sy}(f)$ for $r, s\geq 0$ since $\Omega^{R, r}_{sy}(\Omega^{R, s}_{sy}(f))=\Omega^{R, r+s}_{sy}(f)$. 
This induces an inductive system 
$$\cdots\rightarrow  \HH^{m-p}(A, \Omega^p_{\sy}(A))\rightarrow \HH^{m-p}(A, \Omega^{p+1}_{\sy}(A))\rightarrow \cdots\rightarrow \HH^{m-p}(A, \Omega^{p+r}_{\sy}(A))\rightarrow \cdots.$$ 
It follows from  \cite[Proposition 3.1]{Wan15a} that if $A$ is a Noetherian algebra over a field $k$ such that the enveloping algebra $A\otimes A^{\op}$ is Noetherian, then the colimit of the above inductive system is isomorphic to the $(m-p)$-th {\it Tate-Hochschild cohomology group}  $$\HH_{\sg}^{m-p}(A, A):= \Hom_{\DD_{\sg}(A\otimes A^{\op})}(A, \Sigma^{m-p}A), \ m-p\in \mathbb Z, $$ where $\DD_{\sg}(A\otimes A^{\op})$ is the singularity category of the enveloping algebra $A\otimes A^{\op}$. Recall that the {\it singularity category} $\DD_{\sg}(A)$ (cf. \cite{Buc, Orl}) of a Noetherian algebra $A$ is defined  as the Verdier quotient of the bounded derived category $\DD^b(\mbox{$A$-$\modu$})$ of finitely generated (left) $A$-modules by the full subcategory $\perf(A)$ consisting of complexes quasi-isomorphic to bounded complexes of finitely generated projective $A$-modules. 

\section{Dg $k[\epsilon_i]/(\epsilon_i^2)$-modules}
\subsection{A construction of dg $k[\epsilon_i]/(\epsilon_i^2)$-modules} \label{section3.1}
Let $A$ be a Noetherian algebra  over a field $k$ such that the enveloping algebra $A\otimes A^{\op}$ is Noetherian. 
For $i\in \mathbb Z$, we denote by   $R_i$  the commutative dg algebra $k[\epsilon_i]/(\epsilon_i^2)$ with trivial differential, where $\epsilon_i$ is of degree $i$. With a slight abuse of notation, we denote by $\epsilon_i$ the kernel of the augmentation $R_i\rightarrow k$. Clearly, $\epsilon_i$ is the one-dimensional graded $k$-vector space concentrated in degree $i$. For a chain complex $X$ of (left) $A$-modules, there is a natural isomorphism of chain complexes between $\Sigma^i X$ and the tensor product $  \epsilon_i\otimes X$ . In what follows,  we will not distinguish between them.

Let $\alpha: X\rightarrow Y$ be a morphism (of degree zero) of chain complexes. Recall that the mapping cone of $\alpha$ is defined as the chain complex   $Cone(\alpha)=\Sigma X\oplus Y$ with differential $\left( \begin{smallmatrix} \Sigma d_X &0 \\ \alpha& d_{Y} \end{smallmatrix}  \right).$ We set $C(\alpha)=Cone(\Sigma^{-1}\alpha)$.  Clearly, $C(\alpha)$ is the chain complex $X\oplus \Sigma^{-1}Y$ with  differential 
 $\left( \begin{smallmatrix} d_X &0 \\ \alpha& \Sigma^{-1}d_{Y}\end{smallmatrix}  \right). $
The complex  $C(\alpha)$ may be depicted as $$\xymatrix{
X  \ar@/^2pc/[rr]^-{\alpha} &  \oplus &   \Sigma^{-1}Y. }$$



The following lemma can be used to construct dg $R_i$-modules. 
\begin{lemma}\label{lemma-new3.1}
Let $\beta: X\rightarrow Y$ be a morphism of dg $A$-modules. Let $\alpha: X\rightarrow \Sigma^{i+1}Y$ be another morphism of dg $A$-modules. Then there is a dg $R_i\otimes A$-module structure on $C(\alpha)$ induced by $\alpha$ and $\beta$. 
\end{lemma}
\begin{proof}
By definition, the complex  $C(\alpha)$ is equal to $(X\oplus \Sigma^iY, \left(\begin{smallmatrix} d& 0\\ \alpha& \Sigma^i d\end{smallmatrix}\right)).$
The graded $R_i$-module structure on $C(\alpha)$ is given as follows: For $x+\Sigma^iy\in X\oplus \Sigma^iY$, the action of  $\lambda+\mu \epsilon_i \in R_i \ (\lambda, \mu\in k)$ is 
$$(\lambda+\mu\epsilon_i)(x+\Sigma^iy)=\lambda x+\Sigma^i (\lambda y+\mu\beta(x)).$$ It is clear that this action is compatible with the differential $\left(\begin{smallmatrix} d& 0\\ \alpha& \Sigma^i d\end{smallmatrix}\right)$. This proves the lemma. \end{proof}

Let $m, p\in \Z_{\geq 0}$ and $f\in \HH^{m-p}(A, \Omega_{\sy}^p(A))$. In Section \ref{section2.2} we have defined two liftings $\vartheta^L(f)$ and $\vartheta^R(f)$ associated to $f$. 
It follows from Lemma \ref{lemma-new3.1} that   $C(\vartheta^L(f))$ and $C(\vartheta^R(f))$ are dg $R_{m-p-1}\otimes A\otimes A^{\op}$-modules. To see this, we take the map $\beta$ in Lemma \ref{lemma-new3.1} to be the natural projection $ \Barr_*(A)\rightarrow \Barr_{\geq p}(A)$ and  $\alpha=\vartheta^L(f)$(resp. $\alpha=\vartheta^R(f)$).    
In particular, as  graded  $R_{m-p-1}\otimes A\otimes A^{\op}$-modules, we have $$ C(\vartheta^L(f))\cong\bigoplus_{i=0}^{p-1}(k\otimes \Barr_i(A))\bigoplus R_{m-p-1}\otimes \Barr_{\geq p}(A)\cong C(\vartheta^R(f)),$$ where $k$ is viewed as the  $R_{m-p-1}$-module concentrated in degree zero and thus $k\otimes \Barr_i(A)$ is an  $R_{m-p-1}\otimes A\otimes A^{\op}$-module concentrated in degree $i$. 
In Section \ref{section2.2} we have also defined two cocycles $\Omega_{\sy}^{L, r}(f)$ and $\Omega_{\sy}^{R, r}(f)$ representing the element $\Omega_{\sy}^r(f)\in\HH^{m-p}(A, \Omega_{\sy}^{p+r}(A))$ for $r\geq 0$.  We note that both $C(\Omega_{\sy}^{L, r}(f))$ and $C(\Omega_{\sy}^{R, r}(f))$ are dg  $R_{m-p-1}\otimes A\otimes A^{\op}$-modules. For this,   we take the map $\beta$ in Lemma \ref{lemma-new3.1} to be the projection $\Barr_*(A)\rightarrow \Omega_{\sy}^{p+r}(A)$ induced by the natural map $\Barr_{p+r}(A)\twoheadrightarrow \Omega_{\sy}^{p+r}(A)$ and $\alpha=\Omega_{\sy}^{L, r}(f)$ (resp. $\alpha=\Omega_{\sy}^{R, r}(f)$). In particular, as  graded  $R_{m-p-1}\otimes A\otimes A^{\op}$-modules, we have $$C(\Omega_{\sy}^{L, r}(f)) \cong \bigoplus_{i\neq p+r}( k\otimes \Barr_i(A))\bigoplus  (\Barr_{p+r}(A)\oplus  \Sigma^{m-p-1}\Omega_{\sy}^{p+r}(A))\cong C(\Omega_{\sy}^{R, r}(f)),$$ where $(\Barr_{p+r}(A)\oplus  \Sigma^{m-p-1}\Omega_{\sy}^{p+r}(A))$ is the graded  $R_{m-p-1}\otimes A\otimes A^{\op}$-module determined by the action $$ \epsilon_{m-p-1}\cdot x:= (-1)^{m-p-1}d_{p+r}(x)\in   \Sigma^{m-p-1}\Omega_{\sy}^{p+r}(A)$$ for any $x\in \Barr_{p+r}(A)$.      
When $r=0$, we get that  $C(f)$ is a dg  $R_{m-p-1}\otimes A\otimes A^{\op}$-module. 

\begin{rem}\label{remark-vartheta-h}
Let $f_1$ and $f_2$ be two different  cocycles  representing $f\in \HH^{m-p}(A, \Omega_{\sy}^p(A))$. Then there exists  $\alpha \in \Hom((\Sigma \overline{A})^{\otimes m-1}, \Omega_{\sy}^p(A))$ such that  $f_1-f_2=\delta(\alpha)$. Define a map 
$$\vartheta^L(\alpha): \Barr_*(A)\rightarrow \Sigma^{m-p-1} \Barr_{\geq p}(A)$$ as follows.  Let $x=a_0\otimes \overline{a}_{1, r}\otimes a_{r+1}\in \Barr_{r}(A)$. If $r<m-1$,   we define $\vartheta^L(\alpha)(x)=0.$ If $r\geq m-1$, we define 
\begin{equation*}
  \begin{split}
    \vartheta^L(\alpha)(x)&=a_0\alpha ^L(\overline{a}_{1, m-1})\otimes \overline{a}_{m, r}\otimes a_{r+1}, \end{split}
\end{equation*}
where  $\alpha^L$ is defined as  in Section \ref{section2.2}. Notice that the identity $f_1-f_2=\delta(\alpha)$ yields   $\vartheta^L(\alpha)d+d\vartheta^L(\alpha)=\vartheta^L(f_1)-\vartheta^L(f_2)$. Thus  
the map  $\left( \begin{smallmatrix}\id  & 0\\  \vartheta^L(\alpha) & \id   \end{smallmatrix} \right): C(\vartheta^L(f_2)) \rightarrow C(\vartheta^L(f_1))$ is an isomorphism of dg  $R_{m-p-1}\otimes A\otimes A^{\op}$-modules with inverse $\left( \begin{smallmatrix}\id & 0\\ -\vartheta^L(\alpha) & \id   \end{smallmatrix} \right): C(\vartheta^L(f_1)) \rightarrow C(\vartheta^L(f_2)).$ This shows that $C(\vartheta^L(f))$ does not depend, up to isomorphism of dg $R_{m-p-1}\otimes A\otimes A^{\op}$-modules,  on the choice of the representatives of $f$. Similar arguments are used to prove that $C(\vartheta^R(f)), C( \Omega_{\sy}^{L, r}(f))$ and $C(\Omega_{\sy}^{R, r}(f))$ are  independent of the choice of the representatives of $f$.  \end{rem}

\begin{lemma}\label{lemma-23}
Let $m\in\Z_{>0}$ and  $p\in \Z_{\geq 0}$. For  $f\in \HH^{m-p}(A, \Omega_{\sy}^p(A))$,  the following assertions hold.
\begin{enumerate}[(i)]
\item $ C(\vartheta^L(f))$ is isomorphic to $C(\vartheta^R(f))$ as dg $R_{m-p-1}\otimes A\otimes A^{\op}$-modules.
\item The morphism of dg $R_{m-p-1}\otimes A\otimes A^{\op}$-modules  $$\widehat{\sigma}_p=\left( \begin{smallmatrix} \id& 0\\ 0 & \sigma_{p}  \end{smallmatrix}  \right): C(\vartheta^L(f))\rightarrow C(f)$$ is an isomorphism in the homotopy category $\KK(R_{m-p-1}\otimes A)$ and in $\KK(R_{m-p-1}\otimes A^{\op}),$ where $\sigma_p: \epsilon_{m-p-1}\otimes \Barr_{p}(A) \rightarrow \epsilon_{m-p-1}\otimes \Omega_{\sy}^p(A) $ is the surjection induced by the augmentation $\tau_p: \Barr_p(A) \twoheadrightarrow \Omega_{\sy}^p(A)$. 
\end{enumerate}
\end{lemma}
\begin{proof} Let us prove assertion (i).  Consider the morphism of chain complexes $\left(\begin{smallmatrix} \id & 0\\ h(f)& \id \end{smallmatrix}\right): C(\vartheta^L(f))\rightarrow C(\vartheta^R(f)),$ where $h(f)$ is the chain homotopy defined in (\ref{equ-chain-hom}). Note that  $ \left(\begin{smallmatrix} \id & 0\\ h(f)  & \id \end{smallmatrix}\right)$  is a morphism of dg  $R_{m-p-1}\otimes A\otimes A^{\op}$-modules since $h(f)$ is a morphism of dg $A\otimes A^{\op}$-modules and is compatible with the action of $\epsilon_{m-p-1}$. In fact,  it is an isomorphism with   inverse $ \left(\begin{smallmatrix} \id &0\\ -h(f)& \id \end{smallmatrix}\right): C(\vartheta^R(f))\rightarrow C(\vartheta^L(f)). $ This proves assertion (i). 

Let us prove assertion (ii).   We claim that $C(\vartheta^L(f))$ is isomorphism to 
$$C(\vartheta^L(0))=\Barr_*(A)\oplus \epsilon_{m-p-1}\otimes \Barr_{\geq p}(A)$$ as dg $R_{m-p-1}\otimes A$-modules.   Indeed, we define a morphism of graded $A$-modules  $$\widetilde{f}: \Barr_*(A)\rightarrow \epsilon_{m-p-1}\otimes \Barr_{\geq p}(A)$$ as $\widetilde{f}(x)=\vartheta^L(f)(a_0\otimes \overline{a_1}\otimes \cdots \otimes \overline{a_{i+1}}\otimes 1)$ for $x=a_0\otimes \overline{a_1}\otimes \cdots \otimes \overline{a_i} \otimes a_{i+1}\in \Barr_i(A)$ and $i>0$.    Notice that we have  $\vartheta^L(f)=d\widetilde{f}-\widetilde{f} d$.  This yields a  morphism of dg $R_{m-p-1}\otimes A$-modules
$$\phi(f)=\left( \begin{smallmatrix} \id& 0\\ \widetilde{f}& \id\end{smallmatrix}  \right):C(\vartheta^L(f))\rightarrow C(\vartheta^L(0)) $$
since $\phi(f)$ is compatible with the action of $\epsilon_{m-p-1}$ and  we have 
$$\left(\begin{smallmatrix} d& 0\\ 0& d\end{smallmatrix}  \right)\left( \begin{smallmatrix} \id& 0\\ \widetilde{f}& \id\end{smallmatrix}  \right)=\left( \begin{smallmatrix} \id& 0\\ \widetilde{f}& \id\end{smallmatrix}  \right)\left( \begin{smallmatrix}d& 0\\ \vartheta^L(f)& d\end{smallmatrix}  \right).$$ It is clear that $\phi(f)$ is an isomorphism with inverse 
$\left( \begin{smallmatrix} \id& 0\\ -\widetilde{f}& \id\end{smallmatrix}  \right): C(\vartheta^L(0))\rightarrow C(\vartheta^L(f)).$
This proves the claim.  Similarly, we have an isomorphism of dg $R_{m-p-1}\otimes A$-modules
$$\psi(f)=\left( \begin{smallmatrix} \id& 0\\ \sigma_p\widetilde{f}& \id\end{smallmatrix}  \right): C(f)\rightarrow C(0)=\Barr_*(A)\oplus\epsilon_{m-p-1}\otimes \Omega_{\sy}^p(A),$$
where $\sigma_p\widetilde{f}$ is the following composition of maps
$$\Barr_*(A)\xrightarrow{\widetilde f} \epsilon_{m-p-1}\otimes \Barr_{\geq p}(A)\xrightarrow{\sigma_p} \epsilon_{m-p-1}\otimes \Omega_{\sy}^p(A).$$ Note that   the following   diagram commutes $$
\xymatrix{
C(\vartheta^L(f))\ar[r]^-{\widehat{\sigma}_p}\ar[d]^-{\cong}_-{\phi(f)} & C(f)\ar[d]_-{\cong}^-{\psi(f)}\\
 C(\vartheta^L(0))\ar[r]^-{\widehat{\sigma}_p}  & C(0).\\
}
$$
To prove that $\widehat{\sigma}_p: C(\vartheta^L(f))\rightarrow C(f)$ is an isomorphism in $\KK(R_{m-p-1}\otimes A)$, it is equivalent to prove that  $\widehat{\sigma}_p: C(\vartheta^L(0))\rightarrow C(0)$ is an isomorphism in $\KK(R_{m-p-1}\otimes A)$. We have a commutative diagram of distinguished triangles in $\KK(R_{m-p-1}\otimes A)$
$$\xymatrix{
\Barr_{<p}(A)\ar[r]\ar[d]^-{=} &  C(\vartheta^L(0))\ar[d]^-{\widehat{\sigma}_p}\ar[r] & \Barr_{\geq p}(A)\oplus \epsilon_{m-p-1}\otimes \Barr_{\geq p}(A)\ar[r]\ar[d]^-{\left(\begin{smallmatrix} \id & 0\\ 0 & \sigma_p \end{smallmatrix}\right)} &   \Sigma \Barr_{<p}(A)\ar[d]^-{=}\\
\Barr_{<p}(A)\ar[r] &  C(0)\ar[r] &   \Barr_{\geq p}(A)\oplus \epsilon_{m-p-1}\otimes \Omega_{\sy}^p(A)\ar[r] &   \Sigma \Barr_{<p}(A).}$$
 It is clear that $\left(\begin{smallmatrix} \id & 0\\ 0 & \sigma_p \end{smallmatrix}\right)$ is an isomorphism in $\KK(R_{m-p-1}\otimes A)$ since we have the following commutative diagram
$$\xymatrix@C=0.001pc{
\Barr_{\geq p}(A)\oplus \epsilon_{m-p-1}\otimes \Barr_{\geq p}(A) \ar[rd]^-{\cong}_-{\left(\begin{smallmatrix} \tau_p & 0\\ 0 & \sigma_p \end{smallmatrix}\right)}   \ar[rr]^-{\left(\begin{smallmatrix} \id & 0\\ 0 & \sigma_p \end{smallmatrix}\right)} & &  \Barr_{\geq p}(A)\oplus \epsilon_{m-p-1}\otimes \Omega_{\sy}^p(A)\ar[ld]_{\cong}^-{\left(\begin{smallmatrix} \tau_p & 0\\ 0 & \id \end{smallmatrix}\right)}\\
& \Omega_{\sy}^p(A) \oplus \epsilon_{m-p-1}\otimes \Omega_{\sy}^p(A) }$$
where $\left(\begin{smallmatrix} \tau_p & 0\\ 0 & \id \end{smallmatrix}\right)$ and $\left(\begin{smallmatrix} \tau_p & 0\\ 0 & \sigma_p \end{smallmatrix}\right)$ are isomorphisms in $\KK(R_{m-p-1}\otimes A)$. This implies that $\widehat{\sigma}_p$ is an isomorphism in $\KK(R_{m-p-1}\otimes A)$.  By a similar argument, we can prove that $\widehat{\sigma}_p$ is an isomorphism in $\KK(R_{m-p-1}\otimes A^{\op})$. This proves assertion (ii). The proof is complete.   \end{proof}

\subsection{Dg modules  arising from  the bullet and circle products}\label{section3.2}
From Section \ref{section2.2}, we have a map   $\Omega_{\sy}^r: \HH^{*}(A, \Omega_{\sy}^p(A))\rightarrow \HH^{*}(A, \Omega_{\sy}^{p+r}(A))$ for $p, r\geq 0$. Recall that the Tate-Hochschild cohomology $\HH_{\sg}^*(A, A)$ is isomorphic to  the colimit of the inductive system 
$$\HH^{*}(A, A)\xrightarrow{\Omega_{\sy}^1} \HH^{*}(A,\Omega_{\sy}^1(A))\xrightarrow{\Omega_{\sy}^1}\cdots \xrightarrow{\Omega_{\sy}^1} \HH^{*}(A, \Omega_{\sy}^p(A))\rightarrow \cdots. $$

Let us recall  the Lie bracket $[\cdot, \cdot]$ on $\HH_{\sg}^*(A, A)$ constructed in \cite{Wan15a, Wan15}. The notations  in the present paper are slightly different from those in \cite{Wan15} since we are using the dg bimodules $\Omega_{\nc}^{R, *}(A)$ instead of $\Omega_{\nc}^{L, *}(A)$.  For $f \in C^{m-p}(A, \Omega_{\sy}^p(A))$ and $g\in C^{n-q}(A, \Omega_{\sy}^q(A))$, set
\begin{equation*}
f\bullet_i g:=\begin{cases} 
 (\id^{\otimes q}\otimes f^R )(\id^{\otimes i-1}\otimes \overline{g}\otimes \id^{\otimes m-i})&\mbox{if} \ 1\leq i\leq m, \\
(\id^{\otimes q+i}\otimes \overline{f}\otimes \id^{\otimes -i})(\id^{\otimes m-1}\otimes g^R) & \mbox{if} \ -q\leq i \leq -1,
\end{cases}
\end{equation*}
where $\overline{f}, f^L$ and $f^R$ are defined as in Section \ref{section2.2},  namely 
$$\overline{f}: (\Sigma \overline{A})^{\otimes m} \xrightarrow{f} \Omega_{\sy}^p(A)\hookrightarrow \Barr_{p-1}(A) \xrightarrow{\pi\otimes \id^{\otimes p-1}\otimes \pi} (\Sigma \overline{A})^{p+1};$$ 
$$f^L: (\Sigma \overline{A})^{\otimes m} \xrightarrow{f} \Omega_{\sy}^p(A)\xrightarrow{(\alpha^L_p)^{-1}} \Omega_{\nc}^{L, p}(A)=A\otimes (\Sigma \overline{A})^{\otimes p}; $$
$$f^R:(\Sigma \overline{A})^{\otimes m} \xrightarrow{f} \Omega_{\sy}^p(A)\xrightarrow{(\alpha^R_p)^{-1}} \Omega_{\nc}^{R, p}(A)= (\Sigma \overline{A})^{\otimes p}\otimes A.$$
Clearly,  we have $f\bullet_i g \in C^{m+n-p-q-1}(A, \Omega_{\nc}^{R, p+q}(A))$.  For instance, $f\bullet_i g\ (i>0)$ is the composition of maps
$$(\Sigma \overline{A})^{\otimes m+n-1}\xrightarrow{\id^{\otimes i-1}\otimes \overline{g}\otimes \id^{\otimes m-i}} (\Sigma \overline{A})^{\otimes m+q}\xrightarrow{\id^{\otimes q}\otimes f^R} (\Sigma\overline A)^{\otimes p+q}\otimes A=\Omega_{\nc}^{R, p+q}(A).$$
Since the isomorphism  $\alpha_p^R: \Omega_{\nc}^{R, p}(A)\xrightarrow{\cong} \Omega_{\sy}^p(A)$ (cf.  Lemma \ref{lemma-forms}) induces an   isomorphism $\alpha_p^R: C^{*}(A, \Omega_{\nc}^{R, p}(A))\xrightarrow{\cong} C^{*}(A, \Omega_{\sy}^p(A))$,  we have 
$\alpha_{p+q}^R(f\bullet_i g) \in C^{m+n-p-q-1}(A, \Omega_{\sy}^{p+q}(A))$.  
  We define 
 $$f\circ g:=\sum_{i=1}^m \alpha_{p+q}^R(f\bullet_i g)-(-1)^{(m-p-1)(n-q-1)}\sum_{i=1}^p \alpha_{p+q}^R(g\bullet_{-i}f); $$
 $$f\bullet g:=\sum_{i=1}^m\alpha_{p+q}^R(f\bullet_i g) +\sum_{i=1}^q\alpha_{p+q}^R(f\bullet_{-i} g); $$ 
 \begin{equation*}
 \begin{split}
 [f, g]:=&f\bullet g-(-1)^{(m-p-1)(n-q-1)} g\bullet f\\
 =&f\circ g-(-1)^{(m-p-1)(n-q-1)} g\circ f. 
 \end{split}
 \end{equation*}

 We remark that when these formulas are applied to elements, additional signs will appear because of the Koszul sign rule. When $p=q=0$, $f\circ g=f\bullet g$ is the usual Gerstenhaber circle product and $[\cdot, \cdot]$ is the usual Gerstenhaber bracket on $C^*(A, A)$.    Then from \cite[Section 4.2]{Wan15}, we get that   $[\cdot, \cdot]$ respects the map $\Omega_{\sy}^r: \HH^{*}(A, \Omega_{\sy}^p(A))\rightarrow \HH^{*}(A, \Omega_{\sy}^{p+r}(A))$. Thus  it induces a  well-defined Lie bracket (still denoted by $[\cdot, \cdot]$)  on $\HH_{\sg}^*(A, A)$.  We have the following very important observation. 
 \begin{lemma}\label{lemma-circ-bullet}
For two cocycles $f\in C^{m-p}(A, \Omega_{\sy}^p(A))$ and $g\in C^{n-q}(A, \Omega_{\sy}^q(A))$,  the following identities hold  in $C^{m+n-p-q-1}(A, \Omega_{\sy}^{p+q}(A))$
$$g\bullet f=\Omega_{\sy}^{R, p}(g) h(f)-h_p^{L, R}(g) \vartheta^R(f),$$
$$g\circ f=\Omega_{\sy}^{R, p}(g)h(f)+h_q^{L, R}(f)\vartheta^R(g),$$
where $h(f)$ is defined in (\ref{equ-chain-hom}) and $h_q^{L, R}(f)$ is defined in Remark \ref{remark-delta-cup}. 
\end{lemma}
\begin{proof}
This follows from the following identities 
$$\Omega_{\sy}^{R, p}(g) h(f)=\sum_{i=1}^n \alpha_{p+q}^R(g\bullet_i f), \quad h_p^{L, R}(g)\vartheta^R(f)=\sum_{i=1}^p \alpha_{p+q}^R(g\bullet_{-i}f).$$
Let us verify these two identities. For this, we have 
\begin{equation*}
\begin{split}
&\Omega_{\sy}^{R, p}(g) h(f)(\overline{a}_{1, m+n-1})\\
=&\sum_{i=1}^n(-1)^{(m-p-1)(i-1)}\Omega_{\sy}^{R, p}(g) (\overline{a}_{1, i-1}\otimes \overline{f}(\overline{a}_{i, i+m-1})\otimes \overline{a}_{i+m, m+n-1})\\
=&\sum_{i=1}^n(-1)^{(m-p-1)(i-1)} \alpha_{p+q}^R(\id^{\otimes p}\otimes g^R)( \overline{a}_{1, i-1}\otimes \overline{f}(\overline{a}_{i, i+m-1})\otimes \overline{a}_{i+m, m+n-1})\\
=&\sum_{i=1}^n\alpha_{p+q}^R(g\bullet_i f)(\overline{\alpha}_{1, m+n-1}).
\end{split}
\end{equation*} 
Similarly, we have 
\begin{equation*}
\begin{split}
&h_p^{L, R}(g)\vartheta^R(f)(\overline{a}_{1,m+n-1})\\
=&(-1)^{(m-p)(n-1)}h_p^{L, R}(g)(\overline{a}_{1,n-1 } \otimes f^R(\overline{a}_{n,  m+n-1}))\\
=&(-1)^{(m-p)(n-1)} \sum_{i=1}^p \alpha_{p+q}^R(\id^{\otimes p-i}\otimes \overline{g}\otimes \id^{\otimes i})(\overline{a}_{1,n-1 } \otimes f^R(\overline{a}_{n,  m+n-1}))\\
=&\sum_{i=1}^p \alpha_{p+q}^R(g\bullet_{-i}f)(\overline{\alpha}_{1, m+n-1}).\end{split}
\end{equation*}
This proves the lemma. 
\end{proof}

 The Yoneda product 
 $$\cup': \HH^{m-p}(A, \Omega_{\sy}^p(A))\otimes \HH^{n-q}(A, \Omega_{\sy}^q(A))\rightarrow \HH^{m+n-p-q}(A, \Omega_{\sy}^{p+q}(A))$$ 
 is given by the composition 
 $$\HH^{m-p}(A, \Omega_{\sy}^p)\otimes \HH^{n-q}(A, \Omega_{\sy}^q)\rightarrow \HH^{m+n-p-q}(A, \Omega_{\sy}^{p}\otimes_A\Omega_{\sy}^q)\rightarrow \HH^{m+n-p-q}(A, \Omega_{\sy}^{p+q}),$$
 where we simply write $\Omega_{\sy}^p$ for $\Omega_{\sy}^p(A)$; and the second morphism is the isomorphism induced by $\mu_{p, q}: \Omega_{\sy}^{p}(A)\otimes_A\Omega_{\sy}^q(A)\xrightarrow{\cong} \Omega_{\sy}^{p+q}(A)$ (cf.  the proof of  Lemma \ref{lemma2.2}).  At the complex level, $\cup'$ is given as follows:  For $f\in C^{m-p}(A, \Omega^p_{\sy}(A))$ and $g\in C^{n-q}(A, \Omega^q_{\sy}(A))$,  
 $$f\cup'g(\overline{a}_{1, m+n})=\mu_{p,q}(f(\overline{a}_{1, m})\otimes_A g(\overline{a}_{m+1, m+n})). $$
 We  defined another cup product $\cup$ in \cite[Section 4]{Wan15}: 
 $$f\cup g=\alpha_{p+q}^R(\id^{\otimes p+q}\otimes \mu)(\id^{\otimes q}\otimes f^R\otimes \id) (\id^{\otimes m}\otimes g^R).$$
 More precisely, $f\cup g$ is the composition of maps
 $$\overline A^{\otimes m+n}\xrightarrow{\id^{\otimes m}\otimes g^R} \overline A^{\otimes m+q}\otimes  A\xrightarrow{\id^{\otimes q}\otimes f^R\otimes \id} \overline A^{\otimes p+q}\otimes A\otimes A\xrightarrow{\id^{\otimes p+q}\otimes \mu}  \Omega^{R, p+q}_{\nc}(A)\xrightarrow{\alpha_{p+q}^R} \Omega_{\sy}^{p+q}(A),$$
 where we simply write $\overline{A}$ for $\Sigma \overline{A}$.  At the cohomology level, the cup product  $\cup'$ is equal to $\cup$ (cf. \cite[Section 4]{Wan15}). We note that $\cup$ is compatible with the map $\Omega_{\sy}^r$. Thus, it induces  a well-defined cup product $\cup'=\cup: \HH_{\sg}^*(A, A)\otimes \HH_{\sg}^*(A, A)\rightarrow \HH_{\sg}^*(A, A)$. 
 
\begin{rem}
It is clear that the  two products $\cup'$ and  $\cup$ at the complex level are not (graded-)commutative. 
But we have  the following identity 
\begin{equation}\label{equ-cup}
\begin{split}
f\cup' g-(-1)^{(m-p)(n-q)} g\cup' f=&(-1)^{m-p}\delta(g\bullet f),\\
f\cup g-(-1)^{(m-p)(n-q)} g\cup f=& (-1)^{m-p} \delta(g\circ f),
\end{split}
\end{equation}
for any $f\in C^{m-p}(A, \Omega_{\sy}^p(A))$ and $g\in C^{n-q}(A, \Omega_{\sy}^q(A))$ such that $\delta(f)=0=\delta(g)$ (cf.    \cite[Proposition 4.4]{Wan15}). This shows that $\cup'=\cup$ is graded-commutative at the cohomology level. 

In the following, we will use the   identities in (\ref{equ-cup})  to construct two  dg $R_{m-p-1}\otimes R_{n-q-1}\otimes A\otimes A^{\op}$-modules $C^L(f, g)$ and $C^R(f, g)$ (see below), which are independent  (up to isomorphism) of the choice of representatives in the cohomology classes of $f$  and $g$ (cf. Lemma \ref{lemma3.6}). We stress that these two dg modules play a crucial role in the proof of Proposition \ref{prop-4.3}, a key step in proving  our main Theorem \ref{thm}.

\end{rem}
  Let $f\in \HH^{m-p}(A, \Omega_{\sy}^p(A))$ and $g\in \HH^{n-q}(A, \Omega_{\sy}^q(A))$, which are represented by the cocycles $f\in C^{m-p}(A, \Omega_{\sy}^p(A))$ and $g\in C^{n-q}(A, \Omega_{\sy}^q(A))$ respectively.   Let us consider the following  three chain complexes associated to $f$ and $g$. For simplicity, we set $r:=m-p-1$ and $s:=n-q-1$. 
The first complex is $C(f, g)$ defined as \begin{equation*}  \xymatrix@C=1pc{ \Barr_*(A)\ar@/^1.8pc/[rr]_-{\vartheta^L(f)}\ar@/^4pc/[rrrr]^-{\vartheta^R(g)} &\oplus  &\epsilon_{r} \otimes \Barr_{\geq p}(A) \ar@/_4.5pc/[rrrr]^-{\id_{\epsilon_r}\otimes \Omega_{\sy}^{R, p}(g)}& \oplus&\epsilon_{s} \otimes \Barr_{\geq q}(A) \ar@/_2.5pc/[rr]^-{\id_{\epsilon_s}\otimes \Omega_{\sy}^{L, q}(f)}& \oplus  &\epsilon_{r+s}\otimes\Omega_{\sy}^{p+q}(A). }
\end{equation*}
The identity $\Omega_{\sy}^{L, q}(f)\vartheta^R(g)=\Omega_{\sy}^{R, p}(g)\vartheta^L(f)$ implies that $C(f, g)$ is a well-defined complex. 
 The second one is $C^{L}(f, g)$ defined as 
 \begin{equation*}
  \xymatrix@C=1pc{  \Barr_*(A)\ar@/^1.8pc/[rr]_-{\vartheta^R(f)}\ar@/^4pc/[rrrr]^-{\vartheta^R(g)}\ar@/_6pc/[rrrrrr]^-{g\bullet f}&\oplus  & \epsilon_{r}\otimes \Barr_{\geq p}(A)\ar@/_4pc/[rrrr]^-{\id_{\epsilon_r}\otimes\Omega_{\sy}^{L, p}(g)}& \oplus&  \epsilon_{s}\otimes \Barr_{\geq q}(A)  \ar@/_2pc/[rr]^-{\id_{\epsilon_s}\otimes\Omega_{\sy}^{L, q}(f)}& \oplus  &\epsilon_{r+s}\otimes\Omega_{\sy}^{p+q}(A).}
\end{equation*}
The first identity in (\ref{equ-cup}) ensures that $C^L(f, g)$ is a complex. 
The third one is $C^R(f, g)$ defined as \begin{equation*}\xymatrix@C=1pc{  \Barr_*(A)\ar@/^1.8pc/[rr]_-{\vartheta^R(f)}\ar@/^4pc/[rrrr]^-{\vartheta^R(g)}\ar@/_6pc/[rrrrrr]^-{g\circ f}&\oplus  & \epsilon_{r}\otimes \Barr_{\geq p}(A)\ar@/_4pc/[rrrr]^-{\id_{\epsilon_r}\otimes\Omega_{\sy}^{R, p}(g)}& \oplus&  \epsilon_{s}\otimes \Barr_{\geq q}(A)  \ar@/_2pc/[rr]^-{\id_{\epsilon_s}\otimes\Omega_{\sy}^{R, q}(f)}& \oplus  &\epsilon_{r+s}\otimes\Omega_{\sy}^{p+q}(A),}
\end{equation*}
The second identity in (\ref{equ-cup}) yields that $C^R(f, g)$ is a complex. 

\begin{lemma}\label{lemma2.1}
 For any $f\in \HH^{m-p}(A, \Omega_{\sy}^p(A))$ and $g\in\HH^{n-q}(A, \Omega_{\sy}^q(A))$,  we have isomorphisms of  complexes  $$C(f, g)\cong C^L(f, g) \cong C^R(f, g).$$  
\end{lemma}
\begin{proof} 
 Let us define a map $s^L(f, g): C^L(f, g)\rightarrow C(f, g)$ as  \begin{equation*}
  s^L(f, g):=\left(\begin{smallmatrix}
    \id &0 &0&0\\
    h(f)&\id&0& 0 \\
    0&0&\id& 0\\
    0&h_{p}^{L, R}(g) &0 & \id
  \end{smallmatrix}\right)
\end{equation*}  
where $h_{p}^{L, R}(g)$ is defined in Remark \ref{remark-delta-cup} and $h(f)$ is defined in (\ref{equ-chain-hom}). 
We have  the following identity 
$$\left(\begin{smallmatrix}
   \id &0 &0&0\\
  h(f)  &\id&0& 0 \\
  0  &0&\id& 0\\
    0&h_{p}^{L, R}(g)&0&  \id
  \end{smallmatrix}\right)\left(\begin{smallmatrix}
    d &0 &0&0\\
    \vartheta^R(f)&d&0& 0 \\
    \vartheta^R(g)&0&d& 0\\
    g\bullet f& \Omega_{\sy}^{L, p}(g)&\Omega_{\sy}^{L, q}(f)&0
  \end{smallmatrix}\right)=\left(\begin{smallmatrix}
    d &0 &0&0\\
   \vartheta^L(f)&d&0& 0 \\
    \vartheta^R(g)&0&d& 0\\
    0& \Omega_{\sy}^{R, p}(g)&\Omega_{\sy}^{L, q}(f)&0
  \end{smallmatrix}\right)\left(\begin{smallmatrix}
   \id &0 &0&0\\
  h(f)  &\id&0& 0 \\
  0  &0&\id& 0\\
    0&h_{p}^{L, R}(g)&0&  \id
  \end{smallmatrix}\right),$$
  since $g\bullet f=\Omega_{\sy}^{R, p}(g) h(f)-h_{p}^{L, R}(g)\vartheta^R(f)$ (cf. Lemma \ref{lemma-circ-bullet}). Here,  we simply write $\id_{\epsilon_r}\otimes \Omega_{\sy}^{R, p}(g)$ (resp. $\id_{\epsilon_s}\otimes \Omega_{\sy}^{R, q}(f)$)   as $\Omega_{\sy}^{R, p}(g)$ (resp. $\Omega_{\sy}^{R, q}(f)$).  
  It follows that $s^L(f, g)$ is a morphism of complexes. 
Note that $s^L(f, g)$ is an isomorphism with inverse 
\begin{equation*}
  s^L(f, g)^{-1}=\left(\begin{smallmatrix}
    \id &0 &0&0\\
    -h(f)&\id&0& 0 \\
    0&0&\id& 0\\
   h_{p}^{L, R}(g)h(f)& -h_{p}^{L, R}(g) &0 & \id
  \end{smallmatrix}\right).
\end{equation*}

Let us prove $C^L(f, g)\cong C^R(f, g)$. Consider the following map 
$s'(f, g): C^L(f, g)\rightarrow C^R(f, g)$ given by  
 \begin{equation*}
  s'(f, g):=\left(\begin{smallmatrix}
    \id &0 &0&0\\
    0&\id&0& 0 \\
    0&0&\id& 0\\
    0&h_p^{L, R}(g)&h_q^{L, R}(f) & \id
  \end{smallmatrix}\right).
\end{equation*}
Note that the following identity holds
$$s'(f, g) \left(\begin{smallmatrix}
    d &0 &0&0\\
    \vartheta^R(f)&d&0& 0 \\
    \vartheta^R(g)&0&d& 0\\
    g\bullet f& \Omega_{\sy}^{L, p}(g)& \Omega_{\sy}^{L, q}(f) & 0
  \end{smallmatrix}\right)=\left(\begin{smallmatrix}
    d &0 &0&0\\
    \vartheta^R(f)&d&0& 0 \\
    \vartheta^R(g)&0&d& 0\\
    g\circ f& \Omega_{\sy}^{R, p}(g)& \Omega_{\sy}^{R, q}(f) & 0
  \end{smallmatrix}\right)s'(f, g),$$ since  by Lemma \ref{lemma-circ-bullet},  we have $g\bullet f+h_p^{L, R}(g) \vartheta^R(f)+h_q^{L, R}(f) \vartheta^R(g)=g \circ f$  and by Remark \ref{remark-delta-cup},  we have  
$$h_q^{L, R}(f) d=\Omega_{\sy}^{L, q}(f)-\Omega_{\sy}^{R, q}(f),\quad h_p^{L, R}(g) d=\Omega_{\sy}^{L, p}(g)-\Omega_{\sy}^{R, p}(g).$$  This implies that $s'(f, g)$  is a morphism of complexes. It is clear $s'(f,g)$ is an isomorphism with inverse $$s'(f, g)^{-1}=\left(\begin{smallmatrix}
    \id &0 &0&0\\
    0&\id&0& 0 \\
    0&0&\id& 0\\
    0& -h_p^{L, R}(g)& -h_q^{L, R}(f) & \id
  \end{smallmatrix}\right).$$    This proves the lemma.
\end{proof}

\begin{rem} It is clear that $C(f, g)$ has a natural dg $R\otimes A\otimes A^{\op}$-module structure, where $R:=R_{m-p-1}\otimes R_{n-q-1}$ is the tensor product of the dg algebras $R_{m-p-1}$ and $R_{n-q-1}.$ 
Then, via the above isomorphisms in Lemma \ref{lemma2.1}, the  complexes $C^L(f, g)$ and $C^R(f, g)$ inherit the structure of a  dg  $R\otimes A\otimes A^{\op}$-module from $C(f, g)$. Hence all the three dg $R\otimes A\otimes A^{\op}$-modules are isomorphic. The tensor product $C(f)\otimes_A C(g)$ is endowed with a natural dg  $R\otimes A\otimes A^{\op}$-module structure. \end{rem}

\begin{lemma}\label{lemma3.6}
 Let $m, n\in \Z_{>0}$. For any $f\in \HH^{m-p}(A, \Omega_{\sy}^p(A))$ and $g\in \HH^{n-q}(A, \Omega_{\sy}^q(A))$, we have a morphism of dg $R\otimes A\otimes A^{\op}$-modules 
  $$\Phi(f, g): C(f, g) \rightarrow C(f)\otimes_AC(g)$$
  such that $\Phi(f, g)$ is an isomorphism in $\KK(R\otimes A)$ and in $\KK(R\otimes A^{\op})$. 
\end{lemma}
\begin{proof} Set $r:=m-p-1$ and $s:=n-q-1$. Let us write down the complex $C(f)\otimes_A C(g)$. Recall that $C(f)$ is the following complex  
\begin{equation*}
\xymatrix{
C(f)=\Barr_*(A)  \ar@/^2pc/[rr]^-{f} &  \oplus &\Sigma^{r}\Omega_{\sy}^p(A).
}
\end{equation*}
 Then $C(f)\otimes_A C(g)$ is depicted by  the following diagram
\begin{equation*}
\xymatrix@C=1pc{B_*\otimes_A B_*\ar@/^2pc/[rr]_-{f\otimes_A\id}\ar@/^4pc/[rrrr]_-{\id\otimes_A g} &\oplus& \Sigma^{r}\Omega^p\otimes_A B_*\ar@/_4pc/[rrrr]^-{\id\otimes_A g}&  \oplus& B_*\otimes_A \Sigma^{s}\Omega^q\ar@/_2pc/[rr]^-{f\otimes_A \id}& \oplus & \Sigma^{r}\Omega^p\otimes_A\Sigma^{s}\Omega^q,
  }
\end{equation*}
where, for simplicity, we write $B_*=\Barr_*(A)$ and $\Omega^p=\Omega_{\sy}^p(A)$. Note that there is a  natural  isomorphism of dg $A\otimes A^{\op}$-modules
$$\widetilde{\mu}_p: \Sigma^r\Omega_{\sy}^p(A)\otimes_A \Barr_*(A)\xrightarrow{\cong} \Sigma^{r} \Barr_{\geq p}(A)$$ defined as the composition of maps
$$\Sigma^{r}\Omega_{\sy}^p(A)\otimes_A \Barr_*(A)\xrightarrow{\Sigma^{r}(\alpha_p^L)^{-1}\otimes_A\id} \Sigma^{r}\Omega_{\nc}^{L, p}(A)\otimes_A \Barr_*(A) \xrightarrow{\cong} \Sigma^{r}\Barr_{\geq p}(A),$$
where $(\alpha_p^L)^{-1}$ is defined in Lemma \ref{lemma-forms} and  the second isomorphism is given by 
$$\Sigma^{r}\Omega_{\nc}^{L, p}(A)\otimes_A (A\otimes (\Sigma \overline{A})^{\otimes i}\otimes A)\xrightarrow{\cong} \Sigma^{r}A\otimes (\Sigma \overline{A})^{\otimes p+i}\otimes A\xrightarrow{\cong} \Sigma^{r}\Barr_{p+i}(A).$$ 
Similarly, we have an isomorphism of dg $A\otimes A^{\op}$-modules
$$\widetilde{\mu}_q: \Barr_*(A)\otimes_A \Sigma^{s}\Omega_{\sy}^q(A)\xrightarrow{\cong} \Sigma^{s} \Barr_{\geq q}(A).$$ Recall that we also have an isomorphism $\mu_{p, q}: \Omega_{\sy}^p(A)\otimes_A \Omega_{\sy}^q(A)\xrightarrow{\cong} \Omega_{\sy}^{p+q}(A)$ from Lemma \ref{lemma2.2}. 
Using the above isomorphisms, we get that  $C(f)\otimes_A C(g)$ is isomorphic to the following complex (denoted by $\widetilde{C}_1(f, g)$)
\begin{equation*}
\xymatrix@C=1pc{\Barr_*(A)\otimes_A \Barr_*(A)\ar@/^2pc/[rr]_-{\widetilde{\mu}_p(f\otimes_A\id)}\ar@/^4pc/[rrrr]_-{\widetilde{\mu}_q(\id\otimes_A g)} &\oplus& \Sigma^{r}\Barr_{\geq p}(A)\ar@/_5pc/[rrrr]_-{\mu_{p, q}(\id\otimes_A g)(\widetilde{\mu}_p)^{-1}}&  \oplus& \Sigma^{s}\Barr_{\geq q}(A)\ar@/_2pc/[rr]_{\mu_{p, q}(f\otimes_A \id)(\widetilde{\mu}_q)^{-1}}& \oplus & \Sigma^{r+s}\Omega_{\sy}^{p+q}(A).
  }
\end{equation*}
Here  the isomorphism $C(f)\otimes_AC(g)\xrightarrow{\cong} \widetilde{C}_1(f, g)$ is given by 
$t_1(f, g):=\left(\begin{smallmatrix}\id & 0 & 0 & 0\\ 0& \widetilde{\mu}_p & 0 & 0\\ 0 & 0 & \widetilde{\mu}_q & 0\\ 0 & 0 & 0 & \mu_{p, q}\end{smallmatrix}\right)$. Via this isomorphism, the complex $\widetilde{C}_1(f, g)$ inherits the structure of a dg $R\otimes A\otimes A^{\op}$-module  from $C(f)\otimes_A C(g)$.  We construct a morphism of graded  $R\otimes A\otimes A^{\op}$-modules 
\begin{equation*}
  t(f, g)=\left(\begin{smallmatrix}
    \Delta_{0, 0} & 0& 0 & 0\\
  0 & \id& 0&0\\
  0&0& \id& 0\\
  0&0&0& \id
  \end{smallmatrix}\right):C(f, g)\rightarrow \widetilde{C}_1(f, g), \end{equation*}
  where $\Delta_{0, 0}$ is defined in Section \ref{section2.1}. We claim that $t(f, g)$ commutes with differentials. Indeed, it is sufficient to verify the following identity 
  $$t(f, g) \left(\begin{smallmatrix}
    d &0 &0&0\\
   \vartheta^L(f)&d&0& 0 \\
    \vartheta^R(g)&0&d& 0\\
    0& \Omega_{\sy}^{R, p}(g)&\Omega_{\sy}^{L, q}(f)&0
  \end{smallmatrix}\right)=\left(\begin{smallmatrix}
    d &0 &0&0\\
   \widetilde{\mu}_p(f\otimes_A\id) &d&0& 0 \\
    \widetilde{\mu}_q(\id\otimes_A g) &0&d& 0\\
    0& \mu_{p, q}(\id\otimes_Ag)(\widetilde{\mu}_p)^{-1}&\mu_{p, q}(f\otimes_A\id)(\widetilde{\mu}_q)^{-1}&0
  \end{smallmatrix}\right) t(f, g). $$
  The above identity holds since we have $$\vartheta^L(f)=\widetilde{\mu}_p(f\otimes_A \id)\Delta_{0, 0},\quad \vartheta^R(g)=\widetilde{\mu}_q(\id\otimes_A g)\Delta_{0, 0};$$
  $$\Omega_{\sy}^{R, p}(g)\widetilde{\mu}_p=\mu_{p, q}(\id\otimes_Ag), \quad \Omega_{\sy}^{L, q}(f)\widetilde{\mu}_q=\mu_{p, q}(f\otimes_A \id).$$ This proves the claim. Therefore, we get a morphism of dg $R\otimes A\otimes A^{\op}$-modules 
   $$\Phi(f, g)=t_1(f,g)^{-1}t(f, g): C(f, g)\rightarrow C(f)\otimes_A C(g).$$

It remains to show that $\Phi(f, g)$ is an isomorphism in $\KK(R\otimes A)$ and in $\KK(R\otimes A^{\op})$. For this, it follows from the proof of Lemma \ref{lemma-23} that  there is an isomorphism of dg $R\otimes A$-modules
$$\left(\begin{smallmatrix} \id & 0& 0 & 0\\ \widetilde{f} & \id& 0 & 0\\
\widetilde{g} & 0 & \id& 0\\
0 & d_{p+q} \Omega_{\sy}^{R, p}(\widetilde{g}) & d_{p+q} \Omega_{\sy}^{L, q}(\widetilde{f})&  \id  \end{smallmatrix} \right): C(f, g) \rightarrow C(0, 0),$$ 
 where $\widetilde{g}(x)=\vartheta^R(g)(x\otimes 1)$ and $\widetilde{f}(x)=\vartheta^L(f)(x\otimes 1)$. Here we leave it to the reader to check that  the above map is indeed a morphism of dg $R\otimes A$-modules. Similarly, we have an isomorphism of dg $R\otimes A$-modules 
$C(f)\otimes_A C(g)\xrightarrow{\cong}  C(0)\otimes_A C(0).$  Moreover,  the following diagram commutes in the category of dg $R\otimes A$-modules
$$\xymatrix{
C(f, g) \ar[r]^-{\Phi(f, g)} \ar[d]^-{\cong} & C(f)\otimes_A C(g)\ar[d]^-{\cong} \\
C(0, 0) \ar[r]^-{\Phi(0, 0)} & C(0)\otimes_A C(0).
}$$
  Thus, to prove that $\Phi(f, g) $ is an isomorphism in $\KK(R\otimes A)$, it is equivalent to prove that $\Phi(0, 0)$ is an isomorphism in $\KK(R\otimes A)$. For this, consider the distinguished triangle  $$ B_{<p}\rightarrow C(0)\rightarrow R_{r}\otimes \Omega_{\sy}^p(A)\rightarrow \Sigma B_{<p}$$ 
in $\KK(R_r\otimes A)$. Applying the tensor functor $-\otimes_{A} C(0)$, we get the triangle 
$$B_{<p}\otimes_A C(0)\rightarrow C(0)\otimes_A C(0)\rightarrow R_{r}\otimes (B_{\geq p}(A)\oplus \Sigma^s \Omega_{\sy}^{p+q}(A))\rightarrow \Sigma B_{<p}$$
in $\KK(R\otimes A)$.  Moreover, we have the following commutative diagram 
$$\xymatrix{
B_{<p}\otimes_A B_*\oplus \Sigma^{s}B_{<p+q} \ar[r] &  C(0)\otimes_A C(0)\ar[r] &  R_{r}\otimes (B_{\geq p}(A)\oplus \Sigma^s \Omega_{\sy}^{p+q}(A))\ar[r] &\\
B_{<p}\oplus \Sigma^{s}B_{<p+q} \ar[u]_-{\cong}^-{\left(\begin{smallmatrix} \Delta & 0\\ 0 & \id\end{smallmatrix}\right)} \ar[r]  & C(0, 0)\ar[u]^-{\Phi(0, 0)}\ar[r] & R_{r}\otimes (B_{\geq p}(A)\oplus \Sigma^s \Omega_{\sy}^{p+q}(A))\ar[r]\ar[u]^-{=}&. }$$ Notice that the morphism $\left(\begin{smallmatrix} \Delta & 0\\ 0 & \id\end{smallmatrix}\right)$ is an isomorphism in $\KK(R\otimes A)$, it follows that $\Phi(0,0)$ is an isomorphism and thus $\Phi(f, g)$ is an isomorphism in $\KK(R\otimes A)$. Similarly, $\Phi(f, g)$ is an isomorphism in $\KK(R\otimes A^{\op})$.    This proves the lemma. 
\end{proof}

\section{$R$-relative derived tensor product}\label{section4}
Let us start with the general setting. Let $k$ be a field.  Let $R$ be a commutative dg $k$-algebra and $E$ be a dg $R$-algebra. The  $R$-relative (unbounded) derived category $\DD_R(E)$ is a $k$-linear category with objects being dg $E$-modules. The morphisms of $\DD_R(E)$ are  obtained from morphisms of dg $E$-modules by the localization with respect to  all $R$-relative quasi-isomorphisms, i.e. all morphisms $s:L\rightarrow M$ of dg $E$-modules whose restriction to $R$ is a homotopy equivalence.  For instance, the $k$-relative derived category $\DD_k(E)$ of the dg $k$-algebra $E$ coincides with the usual derived category $\DD(E)$. The  $R$-relative derived category $\DD_R(R )$ is  the homotopy category $\KK(R )$ of $R$.  We also consider the  $R$-relative bounded derived category $\DD_R^b(E)$, which is by definition  the full subcategory of $\DD_R(E)$ consisting of those objects $X$ such that there are only finitely many integers $i$ such that $H_i(X)\neq 0$. For more details on $R$-relative derived categories, we refer to \cite{Kel2, Kel}. 

Let $A$ and $B$ be two associative $k$-algebras. Let $X$ be a dg  $R\otimes A\otimes B^{\op}$-module. Then we have the $R$-relative derived
tensor product induced by $X$, in the sense of \cite{Deli},
$X\otimes_{R\otimes B}^{\bbL, R}- : \DD_{R}( R\otimes B)\rightarrow \DD_{R}( R\otimes A).$
\begin{rem}\label{rem-relative}
From \cite[Section 7]{Kel2}, it follows that
$X\otimes_{R\otimes B}^{\bbL, R}-\cong \mathbf{p}_{rel} X\otimes_{R\otimes B}-,$
where the dg  $R\otimes A\otimes B^{\op}$-module $\mathbf{p}_{rel} X$ is  $R$-relatively quasi-isomorphic to $X$ and {\it $R$-relatively closed} as a  dg $R\otimes B^{\op}$-module, i.e. $\Hom_{\KK(R\otimes B^{\op})}(\mathbf{p}_{rel}X, M)\cong \Hom_{\DD_{R}(R\otimes B^{\op})}(\mathbf{p}_{rel}X, M)$
for any dg $R\otimes B^{\op}$-module $M$.  For instance, we have the isomorphism $\Hom_{\KK(R\otimes B^{\op})}(R\otimes B, M)\cong \Hom_{\DD_R(R\otimes B^{\op})}(R\otimes B, M)$ for any dg $R\otimes B^{\op}$-module $M$, and hence $(R\otimes B)\otimes_{R\otimes B}^{\bbL, R}-\cong (R\otimes B )\otimes_{R\otimes B}-$ (cf. \cite[Section 7.4]{Kel2}). 
\end{rem}
\begin{lemma}\label{lemma-basic}
Let $X$ be a dg $R$-module and $P$ be a ($k$-relatively) closed dg $B^{\op}$-module. Then $X\otimes P$ is $R$-relatively closed as a  dg $R\otimes B^{\op}$-module, namely, we have \begin{equation}\label{equ-lemma}\Hom_{\KK(R\otimes B^{\op})}(X\otimes P, M)\cong \Hom_{\DD_{R}(R\otimes B^{\op})}(X\otimes P, M),\end{equation}  for any  dg $R\otimes B^{\op}$-module $M$. As a consequence,  $(X\otimes P)\otimes_{B\otimes_k R}^{\bbL, R}-\cong (X\otimes P)\otimes_{B\otimes_k R}-.$
\end{lemma}
\begin{proof} It suffices to show that $\Hom_{\KK(R\otimes B^{\op})}(X\otimes P, M)=0$ when $M$ is an $R^{\op}$-contractible (i.e. $M\cong 0$ in $\KK(R^{\op})$) dg $R\otimes B^{\op}$-module. For this, let us write $\mathcal Hom$ for the Hom-complexes.  Equivalently, we need to  show that the complex  $\mathcal Hom_{R\otimes B^{\op}}(X\otimes P, M)$ is acyclic when $M$ is  $R$-contractible. Since we  have 
$$\mathcal Hom_{R\otimes B^{\op}}(X\otimes P, M)\cong\mathcal Hom_{R}(X, \mathcal Hom_{B^{\op}}(P, M)), $$ 
to prove that $\mathcal Hom_{R\otimes B^{\op}}(X\otimes P, M)$ is acyclic, it suffices to prove that  $ \mathcal Hom_{B^{\op}}(P, M)$ is $R$-contractible. Clearly, this holds for $P=B$ and is inherited by shifts and arbitrary coproducts  (because products of $R$-contractible dg $R$-modules are still $R$-contractible). This is also inherited by extensions that split in the category of graded $B^{\op}$-modules. Therefore, it holds for any closed dg $B^{\op}$-modules. This proves the lemma. 
\end{proof}

\begin{rem}
We would  like to thank the referee for providing a shorter proof of Lemma \ref{lemma-basic} and thank Keller for pointing out that this lemma holds for any closed dg $B^{\op}$-module $P$. 
\end{rem}

\begin{prop}\label{lemma-relatively-closed}
Let $m>0$. For  a cocycle $f\in C^{m-p}(A, \Omega_{\sy}^p(A))$ (i.e. $\delta(f)=0$),  all  the three dg modules   $C(f), C(\vartheta^L(f))$ and $C(\vartheta^R(f))$ are $R_{m-p-1}$-relatively closed as  dg $R_{m-p-1}\otimes A$-modules and as dg  $R_{m-p-1}\otimes A^{\op}$-modules. 
\end{prop}
\begin{proof}
It follows from Lemma \ref{lemma-23} that  all the three  modules are isomorphic in $\KK(R_{m-p-1}\otimes A)$ and in $\KK(R_{m-p-1}\otimes A^{\op})$.  Therefore, it is sufficient to prove that $C(\vartheta^L(f))$ is $R_{m-p-1}$-relatively closed as a dg $R_{m-p-1}\otimes A$-module and as a dg $R_{m-p-1}\otimes A^{\op}$-module. From the proof of Lemma \ref{lemma-23}, we have an isomorphism of dg $R_{m-p-1}\otimes A$-modules
$$\phi(f): C(\vartheta^L(f))\xrightarrow{\cong} C(\vartheta^L(0))=\Barr_*(A)\oplus\epsilon_{m-p-1}\otimes \Barr_{\geq p}(A).$$ Let us prove that $C(\vartheta^L(0))$ is $R_{m-p-1}$-relatively closed in $\KK(R_{m-p-1}\otimes A)$. For this, we have a distinguished triangle in $\KK(R_{m-p-1}\otimes A)$
$$\Barr_{0, p-1}(A)\rightarrow C(\vartheta^L(0))\rightarrow R_{m-p-1}\otimes \Barr_{\geq p}(A)\rightarrow \Sigma\Barr_{0, p-1}(A).$$
By Lemma \ref{lemma-basic}  both $\Barr_{0, p-1}(A)$ and $R_{m-p-1}\otimes \Barr_{\geq p}(A)$ are $R_{m-p-1}$-relatively closed in $\KK(R_{m-p-1}\otimes A)$. Hence so is $C(\vartheta^L(0))$. This proves that $C(\vartheta^L(f))$ is $R_{m-p-1}$-relatively closed in $\KK(R_{m-p-1}\otimes A)$.  Similarly, we can prove that it is also $R_{m-p-1}$-relatively closed in $\KK(R_{m-p-1}\otimes A^{\op})$. This prove the proposition.   
\end{proof}




\begin{prop}\label{lemma-plus}
Let $m>0$ and $p\geq 0$. Let $f, g\in \HH^{m-p}(A, \Omega_{\sy}^p(A))$. Then we have a morphism of dg $R_{m-p-1}\otimes A\otimes A^{\op}$-modules $$\Psi(f, g): C(f)\otimes_{R_{m-p-1}\otimes A} C(g)\xrightarrow{\cong}
 C(\Omega_{\sy}^{p}(f+g))$$
 such that $\Psi(f, g)$ is an isomorphism in $\KK(R_{m-p-1}\otimes A)$ and in $\KK(R_{m-p-1}\otimes A^{\op})$. 
Here the map $\Omega_{\sy}^p: \HH^{m-p}(A, \Omega_{\sy}^p(A))\rightarrow \HH^{m-p}(A, \Omega_{\sy}^{2p}(A))$ is  defined in  Section \ref{section2.2}. 
\end{prop}
\begin{proof}  
For simplicity, we write $\Barr_*(A)$ as $B_*$ throughout this proof.  
Then the dg module  $C(f)\otimes_{R_{m-p-1}\otimes A} C(g)$ is illustrated  as follows \begin{equation*}\xymatrix@C=1pc{
 (B_*\ar@/^ 2pc/[rr]^-{f}&  \oplus  &   \Sigma^{m-p-1}\Omega_{\sy}^p(A) ) & \bigotimes_{R_{m-p-1}\otimes A}&   (B_*\ar@/^ 2pc/[rr]^-{g}&   \oplus    &\Sigma^{m-p-1}\Omega_{\sy}^p(A)).}
 \end{equation*} 
 We claim that $C(f)\otimes_{R_{m-p-1}\otimes A} C(g)$  is isomorphic to the following dg  $R_{m-p-1}\otimes A\otimes A^{\op}$-module $C_1(f, g)$
 \begin{equation*}
  \xymatrix@R=0.0000001pc{& \\
   B_*\otimes_A B_*\ar@/^2pc /[rr]^{f\otimes_A d_p+d_p\otimes_A g }& \oplus   & \Sigma^{m-p-1} \Omega_{\sy}^p(A)\otimes_A\Omega^p_{\sy}(A).}
 \end{equation*}
 Indeed,  as graded $R_{m-p-1}\otimes A\otimes A^{\op}$-modules, we have  $$C(f)\cong\bigoplus_{i\neq p}(k\otimes B_i)\bigoplus  (B_{p}\oplus \Sigma^{m-p-1}\Omega_{\sy}^p(A)).$$  Thus,  as  graded $R_{m-p-1}\otimes A\otimes A^{\op}$-modules, we have
\begin{eqnarray*}
\lefteqn{
  C(f)\otimes_{R_{m-p-1}\otimes A} C(g)}\\
    &\cong&\left(\bigoplus_{i\neq p}(k\otimes B_i)\bigoplus  (B_{p}\oplus \Sigma^{m-p-1}\Omega_{\sy}^p(A))\right)\bigotimes_{R_{m-p-1}\otimes A}\\
    && \left(\bigoplus_{i\neq p}(k\otimes B_i)\bigoplus  (B_p\oplus \Sigma^{m-p-1}\Omega_{\sy}^p(A))\right)\\
    &\cong &\left(\bigoplus_{i\neq p}B_i\otimes_{A}\bigoplus_{i\neq p}B_i\right)\bigoplus
   \left( \bigoplus_{i\neq p}B_i\otimes_A B_{ p}\right)\bigoplus\\
   &&\left(B_{ p}\otimes_A\bigoplus_{i\neq p}B_i \right)\bigoplus
   \left(B_p\otimes_A B_p\oplus \Sigma^{m-p-1}\Omega_{\sy}^p(A)\otimes_A \Omega_{\sy}^p(A)   \right)\\
   &\cong& B_*\otimes_A B_*\bigoplus  \Sigma^{m-p-1}\Omega_{\sy}^p(A)\otimes_A \Omega_{\sy}^p(A) \\&\cong& C_1(f, g),
\end{eqnarray*}
where the second isomorphism comes from the following isomorphisms 
\begin{equation}\label{equation-epsilon}
\begin{split}
(k\otimes B_i)\otimes_{R_{m-p-1}\otimes A}(k\otimes B_j) \cong& B_i\otimes_AB_j;\\
X_p\otimes_{R_{m-p-1}\otimes_A} X_p\cong& B_p\otimes_A B_p\oplus \Sigma^{m-p-1}\Omega^p_{\sy}(A)\otimes_A \Omega^p_{\sy}(A);\\
X_p \otimes_{R_{m-p-1}\otimes A} (k\otimes B_i)\cong& B_p\otimes_A B_i, \quad (k\otimes B_i)\otimes_{R_{m-p-1}\otimes A} X_p\cong B_i\otimes_A B_p;
\end{split}
\end{equation}
where $X_p:=B_p\oplus \Sigma^{m-p-1}\Omega_{\sy}^p(A).$ 
The first isomorphism in (\ref{equation-epsilon}) is due to the fact that $k\otimes_{R_{m-p-1}}k\cong k$. Let us prove the second isomorphism in (\ref{equation-epsilon}).  For this, we have a short exact sequence of dg $R_{m-p-1}\otimes A\otimes A^{\op}$-modules
\begin{equation}\label{equation-epsilon1} 0\rightarrow \epsilon_{m-p-1}\otimes \Sigma^{-1}\Omega_{\sy}^{p+1}(A) \rightarrow R_{m-p-1}\otimes B_p \rightarrow X_p\rightarrow 0.
\end{equation}
Applying the tensor functor $-\otimes_{R_{m-p-1}\otimes A} X_p$ to (\ref{equation-epsilon1}), we obtain the following exact sequence 
$$ 0 \rightarrow \Sigma^{-1}\Omega_{\sy}^{p+1}(A)\otimes_A \epsilon_{m-p-1} X_p\rightarrow B_p\otimes_A X_p\rightarrow X_p\otimes_{R_{m-p-1}\otimes A}X_p\rightarrow 0.$$
This implies that 
\begin{equation*}
\begin{split}
X_p\otimes_{R_{m-p-1}\otimes A}X_p\cong& (B_p\otimes_A X_p)/ (\Sigma^{-1}\Omega_{\sy}^{p+1}(A)\otimes_A \epsilon_{m-p-1} X_p)\\
\cong&  (B_p\otimes_A X_p)/((\Sigma^{-1}\Omega_{\sy}^{p+1}(A)\otimes_A \Sigma^{m-p-1} \Omega_{\sy}^p(A))\\
\cong& B_p\otimes_AB_p\oplus \Sigma^{m-p-1}\Omega_{\sy}^p(A)\otimes_A \Omega_{\sy}^p(A).
\end{split}
\end{equation*} 
This shows  the second isomorphism in (\ref{equation-epsilon}).  Similarly, applying the functor $-\otimes_{R_{m-p-1}\otimes A} (k\otimes B_i)$ (resp. $(k\otimes B_i)\otimes_{R_{m-p-1}\otimes A}-$) to (\ref{equation-epsilon1}), we may get the third isomorphism (resp. the forth isomorphism) in (\ref{equation-epsilon}).  Then, from the construction of the tensor product of dg modules, we see that  the differential is exactly given by that  of $C_1(f, g)$. This proves  the claim. 

We identity $C(f)\otimes_{R_{m-p-1}\otimes A} C(g)$ with $C_1(f, g)$ via the above isomorphism. Let us denote by $C_1'(f, g)$ the following dg $R_{m-p-1}\otimes A\otimes A^{\op}$-module \begin{equation*}
 \xymatrix@C=0.001pc{
 B_*\ar@/^ 2pc/[rr]^-{\mu_{p, q}(f\otimes_Ad_p+d_p\otimes_A g) \Delta_{0, 0} }&  \oplus &
 \Sigma^{m-p-1}\Omega_{\sy}^{2p}(A),}
 \end{equation*} where 
   the coproduct $\Delta_{0, 0}: B_*\rightarrow B_*\otimes_A B_*$ and the isomorphism $\mu_{p, q}: \Omega_{\sy}^p(A)\otimes_A \Omega_{\sy}^q(A)\xrightarrow{\cong} \Omega_{\sy}
^{p+q}(A)$ are defined  in  Section  \ref{section2.1}.  
 Note that we have  a morphism of dg $R_{m-p-1}\otimes A\otimes A^{\op}$-modules  $$\Psi'=\left(\begin{smallmatrix} \Delta_{0, 0} & 0\\ 0& \mu^{-1}_{p, q}\end{smallmatrix}\right):  C'_1(f, g)\rightarrow C_1(f, g)= C(f)\otimes_{R_{m-p-1}\otimes A} C(g).$$ 
By Remark \ref{remark-delta-cup}, we have  that $$\mu_{p, p}(f\otimes_A d_p+d_p \otimes_A g)\Delta_{0, 0}=\Omega_{\sy}^{L, p}(f)+\Omega_{\sy}^{R, p}(g)=\Omega_{\sy}^{L, p}(f)+\Omega_{\sy}^{L, p}(g)=\Omega_{\sy}^{p}(f+g)$$ in $\HH^{m-p}(A, \Omega_{\sy}^{2p}(A))$.    
Thus  $C_1'(f, g)$ is isomorphic to $C(\Omega_{\sy}^p(f+g))$ since $C(\Omega_{\sy}^p(f+g))$ does not  depend on the choice of the representatives (cf. Remark \ref{remark-vartheta-h}). This yields a morphism of dg $R_{m-p-1}\otimes A\otimes A^{\op}$-modules 
$$\Psi(f, g): C(\Omega_{\sy}^p(f+g))\rightarrow C(f)\otimes_{R_{m-p-1}\otimes A} C(g). $$ It remains to prove that $\Psi$ is an isomorphism in $\KK(R_{m-p-1}\otimes A)$ and in $\KK(R_{m-p-1}\otimes A^{\op})$.  From the proof of Lemma \ref{lemma-23}, it follows that we have the following commutative diagram
$$\xymatrix{
C(\Omega_{\sy}^{p}(f+g)) \ar[d]^-{\cong}_-{\psi(\Omega_{\sy}^p(g))} \ar[r]^-{\Psi(f, g)} & C(f)\otimes_{R_{m-p-1}\otimes A} C(g)\ar[d]_-{\cong}^-{\id\otimes_{R_{m-p-1}\otimes A} \psi(g)} \\
C(\Omega^p_{\sy}(f)) \ar[r]^-{\Psi(f, 0)}\ar[d]^-{\cong}_-{\psi(\Omega_{\sy}^p(f))}& C(f)\otimes_{R_{m-p-1}\otimes A} C(0)\ar[r]^-{=} & C_1(f, 0) \ar[d]_-{\cong} ^-{\left(\begin{smallmatrix} \id & 0 \\ \psi(f)\otimes \id & \id \end{smallmatrix}\right)}\\
C(0) \ar[r]^-{\Psi(0, 0)} & C(0)\otimes_{R_{m-p-1}\otimes A} C(0)\ar[r]^-{=} & C_1(0, 0)
}$$
in $\KK(R_{m-p-1}\otimes A)$. Thus, to prove that $\Psi(f, g)$ is an isomorphism in $\KK(R_{m-p-1}\otimes A)$, it is equivalent to prove that $\Psi(0, 0)$ is an isomorphism. Note that the latter follows from the following commutative diagram of distinguished triangles in $\KK(R_{m-p-1}\otimes A)$
$$ \xymatrix{
B_{<2p} \ar[r] \ar[d]_{\cong}^-{\Delta}& C(0) \ar[r]\ar[d]^-{\Psi(0, 0)} & \Barr_{\geq 2p}\oplus \Sigma^{m-p-1}\Omega_{\sy}^{2p}(A) \ar[r] \ar[d]_-{\cong}^-{\left(\begin{smallmatrix} \Delta & 0\\ 0 & \id\end{smallmatrix}\right)}& \Sigma B_{<2p}\ar[d]_-{\cong}^-{\Sigma \Delta}\\
(B_*\otimes B_*)_{< 2p}\ar[r] & C_1(0, 0) \ar[r] & (B_*\otimes B_*)_{\geq 2p}\oplus \Sigma^{m-p-1}\Omega_{\sy}^{2p}(A) \ar[r] &\Sigma(B_*\otimes B_*)_{< 2p}. }$$ Similarly, we can prove that $\Psi(f, g) $ is an isomorphism in $\KK(R_{m-p-1}\otimes A^{\op})$. 
This proves the proposition.
 \end{proof}

 \section{Singular infinitesimal deformation theory} In this section, we follow \cite{Kel} and  develop the  singular infinitesimal deformation theory of the identity bimodule. Let $k$ be a field and $R$ be an augmented commutative dg $k$-algebra. We denote by $\bfn$ the kernel of the augmentation $R\rightarrow k$. We always suppose that $\dim_k\bfn<\infty$. 

Let $A$ be a Noetherian $k$-algebra such that the enveloping algebra $A\otimes A^{\op}$ is Noetherian. Let $\DD^{b, \mathrm{Proj}}(A\otimes A^{\op})$ be the full subcategory of $\DD^b(A\otimes A^{\op})$ formed by all the complexes   quasi-isomorphic to bounded complexes $X$  of (not necessarily finitely generated) $A$-$A$-bimodules such that  each component $X_i$ of $X$ is projective as a left $A$-module and as a right $A$-module.  For instance, $A\in \DD^{b, \mathrm{Proj}}(A\otimes A^{\op})$. 

Let $\DD_R^-(R\otimes A\otimes A^{\op})$ be the $R$-relative right bounded derived category of dg $R\otimes A\otimes A^{\op}$-modules (cf. Section \ref{section4}). We consider its full subcategory   $\DD_{R, cl}^{-}(R\otimes A\otimes A^{\op})$ formed by all the objects $X\in \DD_R^-(R\otimes A\otimes A^{\op})$ satisfying the following two conditions
\begin{enumerate}[(i)]
\item $X$ is $R$-relatively closed as a dg  $R\otimes A$-module and as a dg $R\otimes A^{\op}$-module,
\item $k\otimes_RX\in \DD^{b, \mathrm{Proj}}(A\otimes A^{\op})$. \end{enumerate}
  Denote by $\PP(R\otimes A\otimes A^{\op})$ the  thick triangulated subcategory of $\DD_{R, cl}^-(R\otimes A\otimes A^{\op})$   generated by all the objects $P$ such that $k\otimes_RP$ is quasi-isomorphic to  bounded complexes of projective $A\otimes A^{\op}$-modules.     We define  the {\it $R$-relative monoidal singularity category} of $A$  as 
$$\DD_{\sg, R}(A\otimes A^{\op}):=\frac{\DD_{R, cl}^-(R\otimes A\otimes A^{\op})}{\PP(R\otimes A\otimes A^{\op})}.$$   

In particular,  the $k$-relative monoidal singularity category  $\DD_{\sg, k}(A\otimes A^{\op})$ coincides with $\frac{\DD_{cl}^-(A\otimes A^{\op})}{\PP(A\otimes A^{\op})}$, where $\DD_{cl}^-(A\otimes A^{\op})$ is the full subcategory of $\DD^{-}(A\otimes A^{\op})$ formed by all the objects $X$ such that each component $X_i$  is projective as a left $A$-module and as a right $A$-module. Thus it is a full triangulated  subcategory of the (partially) completed singularity category $\widehat{\mathrm{Sg}}(A\otimes A^{\op}):=\DD^-(A\otimes A^{\op})/\PP(A\otimes A^{\op})$ defined in \cite[Section 2.1]{Kel18}.  It follows from Lemma 2.2 in loc. cit. that  the singularity category $\DD_{\sg}(A\otimes A^{\op})$ (in the sense of Buchweitz and Orlov) is a full subcategory of $\widehat{\mathrm{Sg}}(A\otimes A^{\op})$ and is also a full subcategory of $\DD_{\sg, k}(A\otimes A^{\op}).$ As a consequence, we have 
$$\HH_{\sg}^i(A, A):=\Hom_{\DD_{\sg}(A\otimes A^{\op})}(A, \Sigma^iA) \cong \Hom_{\DD_{\sg, k}(A\otimes A^{\op})}(A, \Sigma^iA),\quad \mbox{for any $i\in \Z$}.$$


\begin{lemma}
 The $R$-relative monoidal singularity category $\DD_{\sg, R}(A\otimes A^{\op})$ endowed with the $R$-relative tensor product $\otimes_{R\otimes A}^{\mathbb L, R}$  is a monoidal category with the unit object $R\otimes A$.
 \end{lemma}
 \begin{proof}
 Let us first prove that $\DD_{R, cl}^{-}(R\otimes A\otimes A^{\op})$ is a monoidal category.  Since all the objects in $\DD_{R, cl}^{-}(R\otimes A\otimes A^{\op})$ are $R$-relatively closed as dg $R\otimes A$-modules and as dg $R\otimes A^{\op}$-modules, we have $\otimes_{R\otimes A}^{\mathbb L, R}\cong\otimes_{R\otimes A}$.   Let $X$ and $Y$ be two objects in $\DD_{R, cl}^{-}(R\otimes A\otimes A^{\op})$. Then 
 we claim that the object $ X\otimes_{R\otimes A}Y$ is in $\DD_{R, cl}^{-}(R\otimes A\otimes A^{\op})$. Indeed, $X\otimes_{R\otimes A} Y$ satisfies the condition (i). That is, $X\otimes_{R\otimes A}Y$  is $R$-relatively closed as a dg $R\otimes A$-module and as a dg $R\otimes A^{\op}$-module.    The reason is as follows.  Recall from \cite[Section 7]{Kel2} that  $Y$ is $R$-relatively closed as a dg $R\otimes A$-module if and only if $Y$ admits a filtration of dg $R\otimes A$-modules
 $$0=Y_{-1}\subset Y_0\subset Y_1\subset \cdots \subset Y_p \subset \cdots\subset Y$$
 such that 
 \begin{itemize}
 \item $Y$ is the union of the $Y_p, \ p\in \mathbb Z_{\geq 0}$, 
 \item the  inclusion $Y_{p}\subset Y_{p+1}$ splits in the category of graded $R\otimes A$-modules, \ $p\in \mathbb Z_{\geq 0}$,
 \item each quotient $Y_p/Y_{p-1}$ is isomorphic as dg $R\otimes A$-modules to a direct summand of a direct sum of dg modules $M\otimes A$, where $M$ is dg $R$-module. 
 \end{itemize}
 This induces a filtration of $X\otimes_{R\otimes A}Y$ in the category of dg $R\otimes A$-modules
 $$0=X\otimes_{R\otimes A} Y_{-1}\subset X\otimes_{R\otimes A} Y_0\subset \cdots \subset X\otimes_{R\otimes A} Y_p\subset \cdots \subset X\otimes_{R\otimes A}Y.$$
It follows  that  $X\otimes_{R\otimes A}Y$ is $R$-relatively closed as a dg $R\otimes A$-module,  since each quotient $$(X\otimes_{R\otimes A} Y_p)/(X\otimes_{R\otimes A} Y_{p-1})\cong X\otimes_{R\otimes A} Y_p/Y_{p-1}$$ is $R$-relatively closed as a dg $R\otimes A$-module.  The same argument shows that $X\otimes_{R\otimes A}Y$ is also $R$-relatively closed as a dg $R\otimes A^{\op}$-module. This proves that $X\otimes_{R\otimes A}Y$ satisfies the condition (i). It remains to prove that $X\otimes_{R\otimes A}Y$ satisfies  the condition (ii). This follows from the following isomorphisms   \begin{equation}\label{equation-new}
k\otimes_R(X\otimes_{R\otimes A} Y)\cong (k\otimes_R k)\otimes_{R}(X\otimes_{R\otimes A}Y)\cong (k\otimes_R X)\otimes_{ A} (k\otimes_R Y).
\end{equation}
The proof of the claim is complete. Therefore, $\DD^-_{R, cl}(R\otimes A\otimes A^{\op})$ is a monoidal category. The above isomorphisms (\ref{equation-new}) also implies that $\PP(R\otimes A\otimes A^{\op})$ is a tensor ideal, thus we have that $\otimes_{R\otimes A}^{\mathbb L, R}=\otimes_{R\otimes A}$ is well-defined in $\DD_{\sg, R}(A\otimes A^{\op})$.   This proves the lemma.  
 \end{proof}
  
 \begin{rem}\label{rem-last}
Let $f:R\rightarrow S$ be a morphism  of commutative dg $k$-algebras.  Then  we have a well-defined functor $S\otimes_R- : \DD_{\sg, R} (A\otimes A^{\op}) \rightarrow \DD_{\sg, S}(A\otimes A^{\op})$ since $S\otimes_R X\in \PP(S\otimes A\otimes A^{\op})$ for any $X\in \PP(R\otimes A\otimes A^{\op})$.  \end{rem}

\subsection{Singular infinitesimal deformations}\label{section5.1}
Let $k$ be a field and $A$ be a Noetherian  $k$-algebra such that the enveloping algebra $A\otimes A^{\op}$ is Noetherian. Let $R$ be an augmented commutative dg $k$-algebra such that $\dim_k R<\infty$. Let $\bfn$ be the kernel of the augmentation $R\rightarrow k$ with ${\bfn}^2=0.$ For example, $R=R_m (:=k[\epsilon_m]/\epsilon_m^2)$. 

Define a {\it singular infinitesimal deformation} of $A$ as the pair $(L, u)$, where $L$ is an object in  $\DD^-_{R, cl}(A\otimes A^{\op})$ such that the canonical projection $ \bfn\otimes_R L \rightarrow \bfn L$ is an isomorphism in $\DD_{\sg, k}(A\otimes A^{\op})$;  and $u: k\otimes_R L\rightarrow A$ is an isomorphism in $\DD_{\sg, k}(A\otimes A^{\op})$. We also define $\calF$ as the category whose objects are the singular infinitesimal deformations $(L, u)$ of $A$ and morphisms from $(L, u)$ to $(L', u')$ are given by morphisms $v: L\rightarrow L'$ of $\DD_{\sg, R}(A\otimes A^{\op})$ such that $u'\circ(\id_k\otimes_R v)=u$ in $\DD_{\sg, k}(A\otimes A^{\op})$.  That is, the following diagram 
 \begin{equation*}
 \xymatrix@R=2pc{
 k\otimes_RL \ar[r]^-{u} \ar[d]_-{\id_k\otimes_Rv}&  A\\
 k\otimes_RL' \ar[ur]_-{u'}
 }
 \end{equation*}
 commutes in $\DD_{\sg, k}(A\otimes A^{\op})$.  We denote by $\sgDefo(A, R\rightarrow k)$ the set of isomorphism classes of objects of $\calF$ and denote by $\sgDefo'(A, R\rightarrow k)$ the set of isomorphism classes of {\it weak singular deformations} of $A$, i.e. objects  $L$  in $\DD_{\sg, R}(A\otimes A^{\op})$ such that  $\bfn\otimes_R L\xrightarrow{\cong} \bfn L$  and $k\otimes_RL\cong A$ in $\DD_{\sg, k}(A\otimes A^{\op})$.

 Let $(L, u)$ be an object of $\calF$. Since $L$ is a dg $R\otimes A\otimes A^{\op}$-module, we have the exact sequence $0\rightarrow \bfn L \rightarrow L\rightarrow k\otimes_R L\rightarrow 0$ of dg $A\otimes A^{\op}$-modules, which splits as a sequence of dg $k$-modules (since $k$ is a field). Thus,  it gives rise to a distinguished triangle $\bfn L \rightarrow L\rightarrow k\otimes_RL \rightarrow \Sigma \bfn L$ in $\DD^b(A\otimes A^{\op}).$ Since $\bfn \otimes_RL \cong \bfn L$ in $\DD_{\sg, k}(A\otimes A^{\op})$, we have the distinguished triangle $ \xymatrix{ \bfn\otimes_R L  \ar[r] & L\ar[r]  &   k \otimes_RL  \ar[r]^{\epsilon'} &  \Sigma \bfn \otimes_RL }$ in $\DD_{\sg, k}(A\otimes A^{\op})$.  From $\bfn^2=0$, it follows that  $\bfn \otimes_R L \cong \bfn \otimes (k\otimes_R L) $  as dg modules.   We define a morphism $\epsilon(L, u): A\rightarrow \Sigma \bfn\otimes A$ of $\DD_{\sg,k}(A\otimes A^{\op})$ by the following diagram
 \begin{equation*}
 \xymatrix{
 k\otimes_RL  \ar[r]^{\epsilon'}\ar[d]_-{\cong}^{u} & \Sigma \bfn \otimes_RL \ar[r]^-{\cong} &\Sigma\bfn\otimes (k\otimes_R L)\ar[d]_{\cong}^-{\Sigma\bfn\otimes
 u}\\
 A\ar[rr]^-{\epsilon(L, u)} && \Sigma\bfn \otimes A.
 }
 \end{equation*}
 We claim that the morphism $\epsilon(L, u)$ only depends on the isomorphism class of $(L, u)$ in the category $\calF$.
 Indeed,  let $(L', u')\in \calF$ be such that there exists an isomorphism $v: (L, u)\rightarrow (L', u')$ in $\calF$. To simplify the notational burden, we denote by $k\otimes v$ the morphism $\id_k\otimes v$, etc.  Then we have
 the following commutative diagram in $\DD_{\sg,k}(A\otimes A^{\op})$.
 \begin{equation*}
 \xymatrix@R=2pc{
 k\otimes_RL \ar@/_ 4pc/[dd]_-{u}^-{\cong} \ar[r]^{\epsilon'}\ar[d]_-{\cong}^{k\otimes_Rv} & 
\Sigma\bfn \otimes_RL \ar[d]^-{\cong}_-{\Sigma\bfn\otimes_Rv}\ar[r]^-{\cong} &\Sigma\bfn\otimes (k\otimes_R L)\ar[d]^-{\cong}_-{\Sigma\bfn\otimes (k\otimes_Rv)}\ar@/^ 6pc/[dd]^-{\Sigma\bfn \otimes  u}_-{\cong}\\
   k\otimes_RL' \ar[r]^{\epsilon''}\ar[d]^{u'}_-{\cong} & \Sigma \bfn \otimes_RL' \ar[r]^-{\cong} & \Sigma\bfn\otimes (k\otimes_R L')\ar[d]_-{\Sigma\bfn \otimes u'}^-{\cong}\\
 A\ar[rr]^-{\epsilon(L', u')}_-{\epsilon(L, u)} && \Sigma\bfn \otimes A}
 \end{equation*}
 We claim that the morphism  $k \otimes_Rv: k \otimes_R L\rightarrow k \otimes_R L'$ is an isomorphism in $\DD_{\sg, k}(A\otimes A^{\op})$. Indeed, it suffices to prove that $Cone(k\otimes_R v)=0$  in $\DD_{\sg, k}(A\otimes A^{\op})$.  Since  $v:L\rightarrow L'$ is an isomorphism in $\DD_{\sg, R}(A\otimes A^{\op})$, we have that $Cone(v)$ is in $\PP(R\otimes A\otimes A^{\op})$. Thus, $Cone(k\otimes_R v)\cong k\otimes_R Cone(v) \in \PP (A\otimes A^{\op})$.  This proves the claim. The above commutative diagram implies that  $\epsilon(L, u)=\epsilon(L', u')$.   Therefore we obtain a map $\Phi: \sgDefo(A, R\rightarrow k)\rightarrow \Hom_{\DD_{\sg,k}(A\otimes A^{\op})}(A, \Sigma\bfn\otimes A)$ which sends $(L, u)$ to $\epsilon(L, u)$.

We will  construct the map $\Psi:\Hom_{\DD_{\sg, k}(A\otimes A^{\op})}(A, \Sigma \epsilon_m\otimes A)\rightarrow \sgDefo(A, R_m\rightarrow k)$ in the case of $R=R_m$ for $m\in\Z$ as follows. Let $f: A\rightarrow \Sigma\epsilon_m\otimes A  $ be a morphism in $\DD_{\sg,k}(A\otimes A^{\op})$. Take a representative $f'\in \HH^{m+1}(A, \Omega_{\sy}^{p}(A))$. It follows from Proposition \ref{lemma-relatively-closed} that the dg  $R_m\otimes A\otimes A^{\op}$-module  $C(f')$ lies in $\DD_{R_m, cl}^-(R_m\otimes A\otimes A^{\op})$. We claim that the canonical morphism $\epsilon_m\otimes_{R_m} C(f')\rightarrow \epsilon_mC(f')$ is an isomorphism in $\DD_{\sg, k}(A\otimes A^{\op})$. Indeed,  since $\epsilon_m^2=0$, we have that $$\epsilon_m\otimes_{R_m}C(f') \cong \epsilon_m\otimes(k\otimes_{R_m} C(f'))\cong \epsilon_m\otimes \Barr_*(A),$$
where the second isomorphism comes from the fact that $k\otimes_{R_m}M\cong M/\epsilon_m M$ for any dg $R_{m}$-module $M$.  Then the canonical morphism $\epsilon_m\otimes_{R_m} C(f')\rightarrow \epsilon_mC(f')$ is given by the following commutative diagram in $\DD_{\sg,k}(A\otimes A^{\op})$:
\begin{equation*}
  \xymatrix{
 \epsilon_m \otimes_{R_m} C(f') \ar[d]^-{\cong} \ar[r]   & \epsilon_mC(f')\ar[d]^-{\cong}\\
 \epsilon_m\otimes \Barr_*(A) \ar[r]_-{\cong}^-{\rho}& \epsilon_m\otimes \Omega_{\sy}^p(A),
  }
\end{equation*}
where $\rho$ is induced  by the natural projection $d_p: \Barr_p(A)\rightarrow \Omega_{\sy}^p(A)$ which is an isomorphism  in $\DD_{\sg,k}(A\otimes A^{\op})$. Hence  the  canonical morphism $\epsilon_m\otimes_{R_m}C(f')\rightarrow \epsilon_mC(f')$ is an isomorphism in $\DD_{\sg,k}(A\otimes A^{\op})$. This proves the claim. 
Consider the canonical morphism of complexes of $A$-$A$-bimodules $u': k \otimes_{R_m} C(f')\xrightarrow{\cong} \Barr_*(A) \xrightarrow{d_0} A $. Clearly, it is an isomorphism in $\DD_{\sg,k}(A\otimes A^{\op})$. By definition, we obtain that   $(C(f'), u')\in \calF$.  

Let us define $\Psi(f)=(C(f'), u').$ The following claim ensures that $\Psi(f)$ is well-defined. 
\begin{claim}\label{claim-1}
 $\Psi(f)$ is independent of the choice of the representatives of $f\in\HH_{\sg}^{m+1}(A, A)$. 
 \end{claim} 
 \begin{proof}
 Indeed,  let 
 $f''\in C^{m+1}(A, \Omega_{\sy}^q(A))$ be  another representative of $f$.  Without loss of the generality, we may assume $q\geq p$. Since both $f'$ and $f''$ represent the same element $f$, we have $\Omega_{\sy}^{q-p}(f')=f'' $ in $\HH^{m+1}(A, \Omega_{\sy}^q(A))$, where the map  $\Omega_{\sy}^{q-p}: \HH^{*}(A, \Omega_{\sy}^p(A))\rightarrow \HH^{*}(A, \Omega_{\sy}^q(A))$ is  defined in Section \ref{section2.2}. Equivalently, there exists $h\in C^{m}(A, \Omega_{\sy}^q(A))$ such that $\Omega^{q-p}_{\sy}(f')-f''=\delta(h)$. Now we prove that $$(C(f'), u')=(C(f''), u'')\quad \mbox{in $\sgDefo(A, R_m\rightarrow k)$}.$$ For this, let us define a morphism $\rho: C(f') \rightarrow C(f'')$ of $\DD_{\sg, R_m}(A\otimes A^{\op})$ as follows  \begin{equation}\label{definition-rho}\xymatrix@R=3pc@C=3pc{C(f')\ar[rrr]^{\rho}\ar[d]_-{\widehat{\sigma}_p^{-1}}^-{\cong}& &   & C(f'')\\
C(\vartheta^L(f'))\ar[r]^-{\phi} & C(\vartheta^L(\Omega_{\sy}^{L, q-p}(f')))\ar[r]^-{=}  & C(\vartheta^L(f''+\delta(h)))\ar[r]^-{\left(\begin{smallmatrix} \id & 0\\ -\vartheta^L(h)& \id\end{smallmatrix}\right)}_-{\cong}  &  C(\vartheta^L(f'')).\ar[u]_-{\widehat{\sigma_q}}^-{\cong}}
\end{equation}
Let us explain the above morphisms: The morphism $\widehat{\sigma}_p: C(\vartheta^L(f'))\rightarrow C(f')$ defined in Lemma \ref{lemma-23} is an isomorphism in $\DD^b(R_m\otimes A\otimes A^{\op})$, thus it is an isomorphism in $\DD_{\sg, R_m}( A\otimes A^{\op})$;  The morphism $\phi$ is induced by  the natural projection $\Barr_{\geq p}(A)\rightarrow \Barr_{\geq q}(A)$; The morphism  $\left(\begin{smallmatrix} \id & 0\\ -\vartheta^L(h)& \id\end{smallmatrix}\right)$ is an isomorphism with inverse $\left(\begin{smallmatrix} \id & 0\\ \vartheta^L(h)& \id\end{smallmatrix}\right)$, where $\vartheta^L(h)$ is defined in Remark \ref{remark-vartheta-h}.  Let us prove    that $\phi$ is an isomorphism in $\DD_{\sg, R_m}(A\otimes A^{\op})$. For this, it suffices to prove that $Cone(\phi)\in \PP(R_m\otimes A\otimes A^{\op})$, namely both $Cone(\phi)$ and $k\otimes_R Cone(\phi)$ are quasi-isomorphic to    bounded complexes of projective $A\otimes A^{\op}$-modules. 
Since  $\phi$ is surjective, there is a short exact sequence of dg $R_m\otimes A\otimes A^{\op}$-modules
$$0\rightarrow B_{p, q-1}\rightarrow C(\vartheta^L(f'))\xrightarrow{\phi} C(\vartheta^L(\Omega_{\sy}^{L, q-p}(f'))\rightarrow 0,$$ where $B_{p, q-1}$ is the truncated complex $$B_{p, q-1}: \xymatrix{
0\ar[r] & \Barr_{q-1}(A)\ar[r]^-{d_{q-2}} & \Barr_{q-1}(A) \ar[r]^-{d_{q-1}} \ar[r] & \cdots \ar[r]^-{d_{p+1}} & \Barr_p(A) \ar[r] & 0.}$$  
This induces a distinguished triangle in $\DD(A\otimes A^{\op})$
$$ B_{p, q-1}\rightarrow C(\vartheta^L(f'))\xrightarrow{\phi} C(\vartheta^L(\Omega_{\sy}^{L, q-p}(f'))\rightarrow \Sigma B_{p, q-1}.$$ 
It follows that $Cone(\phi)\cong \Sigma B_{p, q-1}$ in $\DD(A\otimes A^{\op})$ and thus $Cone(\phi)$ is quasi-isomorphic to a bounded complex of projective $A\otimes A^{\op}$-modules. Since $$k\otimes_R Cone(\phi)\cong Cone(\phi)/\epsilon_{m-p-1} Cone(\phi)\cong Cone(\id: \Barr_*(A)\rightarrow \Barr_*(A)),$$
we get that $k\otimes_RCone(\phi)\cong 0$ in $\DD(A\otimes A^{\op})$. This prove that $Cone(\phi)\in \PP(R_m\otimes A\otimes A^{\op})$. Thus the morphism $\rho: C(f')\rightarrow C(f'')$ is an isomorphism of $\DD_{\sg, R_m}(A\otimes A^{\op})$.   

Applying the tensor functor $k\otimes_{R_m}-$ to the above diagram (\ref{definition-rho}), we have  the following  commutative diagram in $\DD_{\sg, k}(A\otimes A^{\op})$
\begin{equation*}
\xymatrix{k \otimes_{R_m} C(f') \ar[r]^-{\id \otimes \rho}  \ar[d]^-{u'}& k \otimes_{R_m}C(f'') \ar[dl]^-{u''}\\ A}
\end{equation*}
since $u'$ is the following composition of maps 
 $$u': k\otimes_{R_m}C(f')\rightarrow \Barr_*(A)\xrightarrow{d_0} A. $$
This yields   $(C(f'), u')=(C(f''), u'')$ in $\sgDefo(A, R_m\rightarrow k).$ Therefore,  $\Psi(f)$ is independent of the choice of the representative of $f$. This proves the claim.
\end{proof}  As a consequence, we get a map $\Psi:\Hom_{\DD_{\sg}(A\otimes A^{\op})}(A, \Sigma\epsilon_m\otimes A)\rightarrow \sgDefo(A, R_m\rightarrow k).$

 \begin{prop}\label{prop-iso}
The  map  $\Psi: \Hom_{\DD_{\sg}(A\otimes A^{\op})}(A, \Sigma\epsilon_m\otimes A)\rightarrow \sgDefo(A, R_m\rightarrow k) $ is injective for any $m\in \Z$. 
 \end{prop}
\begin{proof} It is sufficient to prove that $\Phi \Psi=\id$. For this, let $f'\in \HH^{m+1}(A, \Omega_{\sy}^p(A))$ be a representative of $f\in \HH_{\sg}^{m+1}(A, A)$. Then  we have $\Phi\Psi(f)=\Phi(C(f'), u')=f'.$ This proves the proposition. \end{proof}

Note that the group $\Aut_{\DD_{\sg}(A\otimes A^{\op})}(A)$ of automorphisms of $A$ in $\DD_{\sg}(A\otimes A^{\op})$ acts on $\Hom_{\DD_{\sg}(A\otimes A^{\op})}(A, \Sigma\epsilon_{m}\otimes A)$ via $$s\cdot f:= (\Sigma \epsilon_{m}\otimes s^{-1}) f s$$ for $s\in \Aut_{\DD_{\sg}(A\otimes A^{\op})}(A)$ and $f\in \Hom_{\DD_{\sg}(A\otimes A^{\op})}(A, \Sigma\epsilon_m\otimes A).$ The group $\Aut_{\DD_{\sg}(A\otimes A^{\op})}(A)$   acts on $\sgDefo(A, R_m\rightarrow k)$ via $$s\cdot (L, u):=(L, s u).$$ Clearly,  the forgetful map induces a bijection $$\sgDefo(A, R_m\rightarrow k)/\Aut_{\DD_{\sg}(A\otimes A^{\op})}(A)\cong \sgDefo'(A, R_m\rightarrow k). $$  Recall  that $\sgDefo'(A, R_m\rightarrow k)$ is the set of isomorphism classes of weak singular deformations (cf. the second paragraph of Section \ref{section5.1}).  The map $\Psi$ induces an injection
\begin{equation}\label{rem-inja}
\Psi': \Hom_{\DD_{\sg}(A\otimes A^{\op})}(A, \Sigma\epsilon_m\otimes A)/\Aut_{\DD_{\sg}(A\otimes A^{\op})}(A)\hookrightarrow \sgDefo'(A, R_m\rightarrow k).
\end{equation}
\begin{lemma} \label{lemma4.2}
For any $m\in\Z$,  $\Aut_{\DD_{\sg}(A\otimes A^{\op})}(A)$ acts trivially on $\Hom_{\DD_{\sg}(A\otimes A^{\op})}(A, \Sigma \epsilon_{m}\otimes A)$. As a consequence, the following natural map is bijective.
\begin{equation*}
\Hom_{\DD_{\sg}(A\otimes A^{\op})}(A, \Sigma\epsilon_m\otimes A)\rightarrow   \Hom_{\DD_{\sg}(A\otimes A^{\op})}(A, \Sigma\epsilon_m\otimes A)/\Aut_{\DD_{\sg}(A\otimes A^{\op})}(A).\end{equation*}
\end{lemma}
\begin{proof}  For $s\in \Aut_{\DD_{\sg}(A\otimes A^{\op})}(A)$ and $f\in  \Hom_{\DD_{\sg}(A\otimes A^{\op})}(A, \Sigma\epsilon_m\otimes A)$, we need to show that $(\Sigma\epsilon_{m}\otimes s^{-1}) f s=f.$  Since the Yoneda product in $\DD_{\sg}(A\otimes A^{\op})$ corresponds to the cup product in $\HH^*_{\sg}(A, A)$, we have
\begin{equation*}
\begin{split}
(\Sigma\epsilon_{m}\otimes s^{-1}) f s&=s^{-1}\cup f\cup s= f\cup s^{-1}\cup s=f,
\end{split}
\end{equation*}
where the second identity comes from the fact that the cup product is graded commutative. This shows that the action of $\Aut_{\DD_{\sg}(A\otimes A^{\op})}(A)$ is trivial. Hence  the map   $\Hom_{\DD_{\sg}(A\otimes A^{\op})}(A, \Sigma^mA)\rightarrow   \Hom_{\DD_{\sg}(A\otimes A^{\op})}(A, \Sigma^mA)/\Aut_{\DD_{\sg}(A\otimes A^{\op})}(A)$ is bijective.
\end{proof}

\begin{rem}\label{rem-new5.6}
By Proposition \ref{prop-iso} and Lemma \ref{lemma4.2}, we obtain a natural embedding $\Psi': \Hom_{\DD_{\sg}(A\otimes A^{\op})}(A, \Sigma\epsilon_m\otimes A)\hookrightarrow \sgDefo'(A, R_m\rightarrow k)$ for any $m\in \mathbb Z$. We set $$G_A(\epsilon_m)=\Imm\left(\Psi':\Hom_{\DD_{\sg}(A\otimes A^{\op})}(A, \Sigma \epsilon_m\otimes A)\hookrightarrow\sgDefo'(A, R_m\rightarrow k)\right).$$ 
 \end{rem}
\begin{lemma}\label{lemma5.5}
The isomorphism $\Psi': \Hom_{\DD_{\sg}(A\otimes A^{\op})}(A, \Sigma\epsilon_m\otimes A)\xrightarrow{\cong} G_A(\epsilon_m)$ is a monoid isomorphism, where the monoid structure on $\Hom_{\DD_{\sg}(A\otimes A^{\op})}(A, \Sigma\epsilon_m\otimes A)$ is the additive structure;  and the monoid structure on $G_A(\epsilon_m)$ is given by $\otimes_{R_m\otimes A}^{\mathbb L, R_m}$. \end{lemma}
\begin{proof} Let $f, g\in \Hom_{\DD_{\sg}(A\otimes A^{\op})}(A, \Sigma\epsilon_m\otimes A),$ which are represented by two elements $f_1, g_1 \in \HH^{m+1}(A, \Omega_{\sy}^p(A))$ respectively. From
Proposition \ref{lemma-plus}, it follows that
\begin{equation*}
\Psi'(f)\otimes^{\mathbb L, R_m}_{R_m\otimes A} \Psi'(g)\cong
C(f_1)\otimes^{\mathbb L, R_m}_{R_m\otimes A} C(g_1)\cong C(\Omega_{\sy}^{2p}(f_1+g_1)). 
\end{equation*}
Since  $C(\Omega_{\sy}^{2p}(f_1+g_1))=C(f_1+g_1)=\Psi'(f+g)$ in $ \sgDefo'(A, R_m\rightarrow k)$ (cf. Claim \ref{claim-1}), we get that $\Psi'$ is a monoid morphism. This proves the lemma.  \end{proof}

 \subsection{The generalized Lie algebra of a group-valued functor}
 Let $k$ be a field. Let  $A$ be a Noetherian  $k$-algebra such that the enveloping algebra $A\otimes A^{\op}$ is Noetherian. Denote by ${\bf cdg}_k$ the category of finite-dimensional augmented commutative dg $k$-algebras and by ${\bf grp}$ the category of groups.  We define the functor $\sgDPic_A: {\bf cdg}_k\rightarrow {\bf grp} $ sending $R\in {\bf cdg}_k$ to the {\it $R$-relative singular derived Picard group}
\begin{equation*}
\begin{split}
\sgDPic_A (R ):=&\{ L\in \DD_{\sg, R}(A\otimes A^{\op}) \ | \  \mbox{there exists $L'\in\DD_{\sg, R}(A\otimes A^{\op}) $ such that }\\
&\mbox{  $L\otimes_{R\otimes A}^{\bbL, R}L'\cong L'\otimes_{R\otimes A}^{\bbL, R}L\cong R\otimes A$ in
$\DD_{\sg, R}(A\otimes A^{\op})$}\}/\sim.
\end{split}
\end{equation*}
where $\sim$ means isomorphisms in $\DD_{\sg, R}(A\otimes A^{\op})$. A morphism $f:R\rightarrow S$ in ${\bf cdg}_k$ induces the group homomorphism $\sgDPic_A(f): \sgDPic_A (R )\rightarrow \sgDPic_A(S)$ sending $L\in \sgDPic_A( R)$ to $S\otimes_R L\in \sgDPic_A(S)$ (cf. Remark \ref{rem-last}).   Then the {\it generalized Lie algebra}   $\Lie\sgDPic^*_A$ associated to the group-valued functor $\sgDPic_A$ is given by  ($m\in\mathbb Z$)
\begin{equation*}
\begin{split}
\Lie\sgDPic^m_A:=&\Ker(\sgDPic_A(R_m)\rightarrow \sgDPic_A(k))\\
=&\{L \in \sgDPic_A(R_m)  \ | \ \mbox{ $k\otimes_{R_m} L\cong A$ in $\DD_{\sg,k}(A\otimes A^{\op})$}  \}.
\end{split}
\end{equation*}
 
\begin{rem} Recall from Remark \ref{rem-new5.6} that we denote 
$$G_A(\epsilon_m)=\Imm\left(\Psi':\Hom_{\DD_{\sg}(A\otimes A^{\op})}(A, \Sigma \epsilon_m\otimes A)\hookrightarrow\sgDefo'(A, R_m\rightarrow k)\right).$$
By Proposition \ref{prop-iso}, we have $G_A(\epsilon_m)\hookrightarrow \Lie\sgDPic_A^m$ for any $m\in \mathbb Z$  since  $$  C(f)\otimes_{R_m\otimes A}^{\LL, R_m} C(-f)\cong C(-f)\otimes^{\LL, R_m}_{R_m\otimes A} C(f) \cong R_m\otimes A$$ in $\DD_{\sg, R_m}(A\otimes A^{\op})$ and  $k \otimes_{R_m}C(f)\cong A$ in  $\DD_{\sg,k}(A\otimes A^{\op})$.
\end{rem}

It follows from Lemma \ref{lemma5.5}  that $\Psi'$ is a monoid isomorphism. Hence $G_A(\epsilon_m)$ has a $k$-vector space structure inherited from that of $\Hom_{\DD_{\sg}(A\otimes A^{\op})}(A, \Sigma\epsilon_m\otimes A).$  We will define a Lie bracket on $G_A:=\bigoplus_{m\in \Z} G_A(\epsilon_m)$ as follows. Let $L_1$ and $L_2$ represent elements of $G_A(\epsilon_m)$ and $G_A(\epsilon_n)$, respectively. Let $U_i$ be the image of $L_i$ in $\DD_{\sg, R}(A\otimes A^{\op})$ where $R=R_m\otimes R_n$ ($i=1, 2$). Note that $U_i$ are invertible objects of the monoidal category $\DD_{\sg, R}(A\otimes A^{\op}) $ (cf. Lemma \ref{lemma5.5}), namely $U_i\in \sgDPic_A(R)$. Let $V$ be the commutator of $U_1$ with $U_2$, namely $V=U_1U_2U_1^{-1}U_2^{-1}\in \sgDPic_A(R)$. Then Proposition \ref{prop-4.3} below shows that $V$ lies in   $G_A(\epsilon_{m+n})$ under the  morphism $\sgDPic_A(R_{m+n})\rightarrow \sgDPic_A(R)$ induced by the natural embedding $R_{m+n}\hookrightarrow R$. Let us define $[L_1, L_2]:=V\in G_A(\epsilon_{m+n}).$

\begin{prop}\label{prop-4.3}
Let $f\in \HH^{m+1}_{\sg}(A, A)$ and $g\in \HH_{\sg}^{n+1}(A, A)$.
Then the commutator $[\Psi'(f), \Psi'(g)]:=\widehat{\Psi'(f)}\otimes^{\bbL, R_m\otimes R_n}_{R_m\otimes R_n\otimes A_n}\widehat{\Psi'(g)} \otimes^{\bbL, R_m\otimes R_n}_{R_m\otimes R_n\otimes A} \widehat{\Psi'(-f) }\otimes^{\bbL, R_m\otimes R_n}_{R_m\otimes R_n\otimes A}\widehat{\Psi'(-g)}$ equals to $\widehat{\Psi'([f, g])}$ in $G_A(\epsilon_m\otimes\epsilon_n),$ where we write 
\begin{equation*}
\begin{split}
\widehat{\Psi'(f)}:=&R_n\otimes \Psi'(f), \quad  \widehat{\Psi'(g)}:=R_m\otimes \Psi'(g),\\
\widehat{\Psi'([f, g])}:=&(R_m\otimes R_n)\otimes_{R_{m+n}}\Psi'([f, g]).
\end{split}
\end{equation*}
Here $[f, g]$ is the  Lie bracket in $\HH_{\sg}^*(A, A)$ (cf. Section \ref{section3.2}).
\end{prop}
\begin{proof} Note that $\Psi'(f)=C(f)$. Then from Lemma \ref{lemma5.5}, it follows that to verify  the identity $[\Psi'(f), \Psi'(g)]=\widehat{\Psi'([f, g])}$ in $\sgDefo'(A, R_{m+n}\rightarrow k)$ is equivalent to verify the following isomorphism in $\DD_{\sg, R_m\otimes R_n}(A\otimes A^{\op})$ \begin{equation}\label{equ-ghy}
\widehat{\Psi'(f)}\otimes^{\bbL, R_m\otimes R_n}_{R_m\otimes R_n\otimes A}\widehat{\Psi'(g)} \cong \widehat{\Psi'([f, g])}\otimes_{R_m\otimes R_n\otimes A}^{\bbL, R_m\otimes R_n} \widehat{\Psi'(g)}\otimes^{\bbL, R_m\otimes R_n}_{R_m\otimes R_n\otimes A}\widehat{\Psi'(f)}.
\end{equation}
By Lemma \ref{lemma3.6}, the left hand side of (\ref{equ-ghy}) is isomorphic  to $C^R(f, g).$ The right hand side  is 
\begin{equation*}
\begin{split}
\RHS&\cong ((R_m\otimes R_n)\otimes_{R_{m+n}}C([f, g]))\otimes_{R_m\otimes R_n \otimes A}^{\bbL, R_m\otimes
R_n}C^L(g, f)\\
&\cong ((R_m\otimes R_n)\otimes_{R_{m+n}}C([f, g]))\otimes_{A\otimes
R_m\otimes R_n}C^L(g, f)\\
&\cong C([f, g]) \otimes_{R_{m+n}\otimes A} C^L(g, f), 
\end{split}
\end{equation*}
where the first isomorphism follows from Lemma \ref{lemma3.6} and the second one is because of the fact that $((R_m\otimes R_n)\otimes_{R_{m+n}}C([f, g]))$ is $(R_m\otimes R_n)$-relatively closed. Now let us compute $C([f, g]) \otimes_{ R_{m+n}\otimes A} C^L(g, f)$ which is illustrated as follows 
\begin{equation*}
\xymatrix@C=1pc@R=4pc{(B_*\ar@/^ 2pc/[rr]^-{[f, g]} &\oplus&\epsilon_{m+n}\otimes  \Omega_{\sy}^{p+q}(A))& \bigotimes_{ R_{m+n}\otimes A}\\
B_*\ar@/^1.8pc/[rr]_-{\vartheta^R(g)}\ar@/^4pc/[rrrr]_-{\vartheta^R(f)}  \ar@/_6pc/[rrrrrr]^-{f\bullet g}&\oplus  & \epsilon_{n}\otimes B_{\geq q} \ar@/_4pc/[rrrr]^-{\Omega_{\sy}^{L, p+q}(f)}& \oplus&  \epsilon_{m}\otimes B_{\geq p} \ar@/_2pc/[rr]^-{\Omega^{L, p+q}_{\sy}(g)}& \oplus  &\epsilon_{m+n}\otimes \Omega_{\sy}^{p+q}(A),
  }
  \end{equation*}
  where, for simplicity,  we write $\Barr_*(A)$ as $B_*$. 
  As graded modules, we have 
  \begin{equation*}
  \begin{split}
 & C([f, g])\otimes_{ R_{m+n}\otimes A} C^L(g, f)\\
  \cong&  B_*\otimes_A \left (B_*\oplus \Sigma^nB_{\geq q}\oplus \Sigma^m B_{\geq p}\right)\oplus \Sigma^{m+n}\Omega_{\sy}^{p+q}(A)\otimes_A \Omega_{\sy}^{p+q}(A)\\
  \cong & B_*\otimes_A B_*\oplus B_*\otimes_A \Sigma^n B_{\geq q}\oplus B_*\otimes \Sigma^m B_{\geq p}\oplus \Sigma^{m+n}\Omega_{\sy}^{p+q}(A)\otimes_A \Omega_{\sy}^{p+q}(A),   \end{split}
  \end{equation*}
  where in the first identity we use the following two  isomorphisms 
  $$(k\otimes B_i)\otimes_{R_{m+n}\otimes A} M\cong B_i\otimes_A (M/\epsilon_{m+n}M)$$
       $$(B_{p+q}\oplus \Sigma^{m+n} \Omega_{\sy}^{p+q}(A))\otimes_{R_{m+n}\otimes A} M\cong (B_{p+q}\otimes_AM) / (\Sigma^{-1} \Omega_{\sy}^{p+q+1}(A)\otimes_A \epsilon_{m+n}M)$$
for any dg $R_{m+n}\otimes A$-module $M$.  
The proofs of the above two isomorphisms are similar to the ones of the isomorphisms in (\ref{equation-epsilon}).  From the construction of the tensor product of dg modules, we obtain the differential illustrated as follows
   \begin{equation*}
\xymatrix@C=1.3pc@R=4pc{B_*\otimes_A B_*\ar@/^1.8pc/[rr]_-{\id\otimes_A\vartheta^R(g)}\ar@/^4pc/[rrrr]_-{\id\otimes_A\vartheta^R(f)}  \ar@/_7pc/[rrrrrr]^-{[f, g]\otimes_Ad_{p+q}+d_{p+q}\otimes_A(f\bullet g)}&\oplus  & B_*\otimes_A \Sigma^n B_{\geq q}\ar@/_4.5pc/[rrrr]^-{d_{p+q}\otimes_A\Omega_{\sy}^{L, p+q}(f)}& \oplus&  B_*\otimes_A \Sigma^m B_{\geq p} \ar@/_2.5pc/[rr]^-{d_{p+q}\otimes_A\Omega^{L, p+q}_{\sy}(g)}& \oplus  &\Sigma^{m+n}\Omega_{\sy}^{p+q}\otimes_A \Omega_{\sy}^{p+q}.
      }
\end{equation*}
Using the quasi-isomorphism $\Delta=\Delta_{p, q}: B_{\geq p+q}\rightarrow B_{\geq p}\otimes_AB_{\geq q}$ and the isomorphism 
$\mu=\mu_{p+q, p+q}: \Omega_{\sy}^{p+q}(A)\otimes_A \Omega_{\sy}^{p+q}(A)\xrightarrow{\cong} \Omega_{\sy}^{2(p+q)}(A)$, the above dg $R_{m}\otimes R_n\otimes A\otimes A^{\op}$-module is $R_m\otimes R_n$-relatively quasi-isomorphic to the following one (denoted by $C_2(f, g)$)
  \begin{equation*}
  \xymatrix@C=2pc@R=1pc{B_*\ar@/^1.8pc/[rr]_-{\vartheta^R(g)}\ar@/^4pc/[rrrr]_-{\vartheta^R(f)}  \ar@/_6pc/[rrrrrr]^-{H}&\oplus  & \epsilon_{n}\otimes  B_{\geq q}\ar@/_4pc/[rrrr]^-{\widetilde{f}}& \oplus&\epsilon_{m}  \otimes B_{\geq p} \ar@/_2pc/[rr]^-{\widetilde{g}}& \oplus  &\epsilon_{m+n}\otimes \Omega_{\sy}^{2(p+q)}(A)      }
\end{equation*}
where we take
\begin{enumerate}[(i)]
\item $\widetilde{g}=\mu(d_{p+q}\otimes_A \Omega_{\sy}^{L, p+q}(g))\Delta$;
\item $\widetilde{f}=\mu(d_{p+q}\otimes_A \Omega_{\sy}^{L, p+q}(f))\Delta$;
\item $H=\mu([f, g]\otimes_Ad_{p+q}+d_{p+q}\otimes_A(f\bullet g))\Delta.$
\end{enumerate}
From Remark \ref{remark-delta-cup}, we have 
$$H=\Omega_{\sy}^{L, p+q}([f, g])+\Omega_{\sy}^{R, p+q}(f\bullet g).$$
We claim that $C_2(f, g)$ is isomorphic to the following dg $R_{m}\otimes R_n\otimes A\otimes A^{\op}$-module (denoted by $C'_2(f, g)$) 
\begin{equation*}
\xymatrix@C=2pc@R=5pc{B_*\ar@/^1.8pc/[rr]_-{\vartheta^R(g)}\ar@/^4pc/[rrrr]_-{\vartheta^R(f)}  \ar@/_6pc/[rrrrrr]^-{\Omega_{\sy}^{R, p+q}(g\circ f)}&\oplus  & \epsilon_{n}\otimes  B_{\geq q}\ar@/_4pc/[rrrr]^-{\Omega_{\sy}^{R, p+2q}(f)}& \oplus&\epsilon_{m}  \otimes B_{\geq p} \ar@/_2pc/[rr]^-{\Omega_{\sy}^{R, 2p+q}(g)}& \oplus  &\epsilon_{m+n}\otimes \Omega_{\sy}^{2(p+q)}(A).     }
\end{equation*}
Indeed, since both $\widetilde{g}$ and $\Omega_{\sy}^{R, 2p+q}(g)$ represent the cocycle $\Omega_{\sy}^{2p+q}(g)$,  there is a coboundary $g_1: \Barr_{n+2p+q-1}\rightarrow \Omega_{\sy}^{2(p+q)}(A)$ such that $g_1d_{n+2p+q}=\Omega_{\sy}^{R, 2p+q}(g)-\widetilde{g}$.  Note that 
$$\widetilde{g}-\Omega_{\sy}^{R, 2p+q}(g)=\mu(d_{p+q}\otimes_A (\Omega_{\sy}^{L, p+q}(g)-\Omega_{\sy}^{R, p+q}(g)))\Delta,$$  it follows from Remark \ref{remark-delta-cup} that we may choose $g_1=-\mu(d_{p+q}\otimes_A h_{p+q}^{L, R}(g))\Delta$ . Similarly, we take  $$f_1=-\mu (d_{p+q}\otimes_Ah_{p+q}^{L, R}(f))\Delta: \Barr_{m+2q+p-1}(A)\rightarrow \Omega_{\sy}^{2(p+q)}(A).$$ Then we have  $f_1d_{m+2q+p}=\Omega_{\sy}^{R, p+2q}(f)-\widetilde{f}$. Since $[f, g]$ is a cocycle, there is a homotopy $$h=h^{L, R}_{p+q}([f, g]): \Barr_{m+n+p+q-2}(A)\rightarrow \Omega^{2(p+q)}_{\sy}(A)$$ such that $h d_{m+n+p+q-1}=\Omega_{\sy}^{L, p+q}([f, g])-\Omega_{\sy}^{R, p+q}([f, g]).$
Let us construct a morphism of graded $A\otimes A^{\op}$-modules
 $$\rho=\left(\begin{smallmatrix} \id & 0 & 0 & 0\\ 0 &  \id & 0 & 0\\ 0  & 0 & \id  & 0\\h & f_1& g_1& \id \end{smallmatrix}\right): C_2(f, g)\rightarrow C'_2(f, g).$$
The following identity holds 
$$\left(\begin{smallmatrix} \id & 0 & 0 & 0\\ 0 &  \id & 0 & 0\\ 0  & 0 & \id  & 0\\h & f_1& g_1& \id \end{smallmatrix}\right) \left(\begin{smallmatrix}
    d &0 &0&0\\
    \vartheta^R(g)&d&0& 0 \\
    \vartheta^R(f)&0&d& 0\\
    H& \widetilde{f}& \widetilde{g} & 0
  \end{smallmatrix}\right)=\left(\begin{smallmatrix}
    d &0 &0&0\\
    \vartheta^R(g)&d&0& 0 \\
    \vartheta^R(f)&0&d& 0\\
    \Omega_{\sy}^{R, p+q}(g\circ f)& \Omega_{\sy}^{R, p+2q}(f)& \Omega_{\sy}^{R, 2p+q}(g) & 0
  \end{smallmatrix}\right)\left(\begin{smallmatrix} \id & 0 & 0 & 0\\ 0 &  \id & 0 & 0\\ 0  & 0 & \id  & 0\\h & f_1& g_1& \id \end{smallmatrix}\right), $$ since  we have 
  \begin{equation*}
  \begin{split} f_1\vartheta^R(g)+g_1\vartheta^R(f)=&\mu(d_{p+q}\otimes_A h_{p+q}^{L, R}(f)\vartheta^{R}(g))\Delta+\mu(d_{p+q}\otimes_A h_{p+q}^{L, R}(g)\vartheta^{R}(f))\Delta\\
  =&\mu(d_{p+q}\otimes_A (g\circ f-g\bullet f))\Delta\\
  =&\Omega_{\sy}^{R, p+q}(g\circ f-g\bullet f),
  \end{split}
  \end{equation*}
  where the first identity is due to the definition of $\Delta$; the second identity is because of  Lemma \ref{lemma-circ-bullet}; and the third identity follows from Remark \ref{remark-delta-cup}. 
This implies that  $\rho$ is a morphism  of dg $R_m\otimes R_n\otimes A\otimes A^{\op}$-modules. It is clear that $\rho$ is bijective  with inverse 
$\rho^{-1}=\left(\begin{smallmatrix} \id & 0 & 0 & 0\\ 0 &  \id & 0 & 0\\ 0  & 0 & \id  & 0\\ -h & -f_1& -g_1& \id \end{smallmatrix}\right).$ This proves the claim. In conclusion, the right hand side of (\ref{equ-ghy}) is isomorphic to $C_2'(f, g)$ in $\DD_{\sg, R_m\otimes R_n}( A\otimes A^{\op})$.

In order to prove (\ref{equ-ghy}),  now it remains to prove that  $C'_2(f, g)$ is isomorphic  to $C^R(f, g)$ in $\DD_{\sg, R_m\otimes R_n}( A\otimes A^{\op})$.  From Lemma \ref{lemma3.6} and Claim \ref{claim-1},  it follows that   $C^R(f, g)$ as an object of $\DD_{\sg, R_m\otimes R_n}(A\otimes A^{\op})$  does not depend on the choice of the representatives of $f$ and $g$ in $\HH^*_{\sg}(A,A)$. Thus $C^R(f, g)$ is isomorphic to $C^R(\Omega_{\sy}^{R, q}(f), \Omega_{\sy}^{R, p}(g))$ in $\DD_{\sg, R_m\otimes R_n}(A\otimes A^{\op})$. It suffices to show  that $C'_2(f, g)$ is isomorphic to $C^R(\Omega_{\sy}^{R, q}(f), \Omega_{\sy}^{R, p}(g))$. For this, since   $$\Omega_{\sy}^{R, p+q}(g\circ f)=\Omega_{\sy}^{R, p}(g) \circ \Omega_{\sy}^{R, q}(f),\quad \Omega_{\sy}^{R, 2p+q}(g)=\Omega_{\sy}^{R, p+q}(\Omega_{\sy}^{R, p}(g))$$ in $C^*(A, \Omega_{\sy}^{2(p+q)}(A)),$ we have a short exact sequence 
$$0\rightarrow B_{q, p+q-1}\oplus B_{p, p+q-1} \rightarrow C_2'(f, g)\xrightarrow{\pi} C^R(\Omega_{\sy}^{R, q}(f), \Omega_{\sy}^{R, p}(g))\rightarrow 0,$$
where $\pi$ is the natural projection; and $B_{p,q}$ denotes   $$ \xymatrix{
0\ar[r] & \Barr_{q}(A)\ar[r]^-{d_{q}} & \Barr_{q-1}(A) \ar[r]^-{d_{q-1}} \ar[r] & \cdots \ar[r]^-{d_{p+1}} & \Barr_p(A) \ar[r] & 0.}$$
Notice that $k\otimes_{R_m\otimes R_n} C_2'(f, g) \cong \Barr_*(A)\cong k\otimes_{R_m\otimes R_n} C^R(\Omega_{\sy}^{R, q}(f), \Omega_{\sy}^{R, p}(g))$, we obtain that  $$k\otimes_{R_m\otimes R_n} Cone(\pi)\cong Cone(k\otimes_{R_m\otimes R_n} \pi) \cong 0$$ in $\DD(A\otimes A^{\op})$  and thus $Cone(\pi)\in \PP(R_m\otimes R_n\otimes A\otimes A^{\op})$. So we  obtain that   $C'_2(f, g)$ is isomorphic  to $C^R(\Omega_{\sy}^{R, q}(f), \Omega_{\sy}^{R, p}(g))$ in $\DD_{\sg, R_m\otimes R_n}( A\otimes A^{\op})$.  This proves the isomorphism (\ref{equ-ghy}). 
The proof is complete.   
\end{proof}

\begin{cor}\label{cor34}
Let $k$ be a field. Let $A$ be a Noetherian $k$-algebra such that the enveloping algebra $A\otimes A^{\op}$ is Noetherian. Then the isomorphisms $\Psi'_m:   \Hom_{\DD_{\sg}(A\otimes A^{\op})}(A, A\otimes \Sigma\epsilon_m )\rightarrow G_A(\epsilon_m)$ induce an isomorphism of graded Lie algebras between $\HH^{*+1}_{\sg}(A, A)$ and $G_A.$
\end{cor}
\begin{proof} This is a direct consequence of Proposition \ref{prop-4.3}.
\end{proof}

\begin{rem}
 We do not know whether the generalized Lie algebra $\Lie\sgDPic_A^*$ is a graded Lie algebra since if two elements $L_1$  and $L_2$ in $\Lie\sgDPic_A^*$ do not lie in  the subspace $G_A$, then it is not clear whether their commutator lies in $\Lie\sgDPic_A^*$.    But however,  it follows from Corollary \ref{cor34} that the  graded subspace $G_A\subset \Lie\sgDPic_A^*$ is indeed a graded Lie algebra.   Keller in \cite{Kel} proved the identity $[\Psi'(f), \Psi'(g)] =\Psi'([f, g])$ for any  $f, g\in \HH^{*+1}(A, A)$ in a quite different way,  where he used the intrinsic interpretation of the Gerstenhaber bracket by Stasheff \cite{Sta}. 
\end{rem}

\section{The invariance under singular equivalence of Morita type with level}\label{section6} Let $k$ be a field.  Let $A$ and $B$ be two Noetherian $k$-algebras such that the enveloping algebras $A\otimes A^{\op}$ and $B\otimes B^{\op}$ are Noetherian. Let $_AM_B$ and $_BN_A$ be an $A$-$B$-bimodule and a  $B$-$A$-bimodule, respectively. Recall from \cite{Wan14} that $(_AM_B, _BN_A)$ defines a {\it singular equivalence of Morita type with level} $l\in \Z_{\geq 0}$ if the following conditions are satisfied:
\begin{enumerate}
\item $M$ is finitely generated projective as a left $A$-module and as a right $B$-module,
\item $N$ is finitely generated projective as a left $B$-module and as a right $A$-module,
\item there exist isomorphisms $M\otimes_B N\cong\Omega_{\sy}^l(A)$ in $(A\otimes A^{\op})$-$\underline{\modu}$,
and $N\otimes_A M\cong \Omega_{\sy}^l(B)$ in $(B\otimes B^{\op})$-$\underline{\modu}$, where $A\otimes A^{\op}$-$\underline{\modu}$ denotes the stable module category of $A$-$A$-bimodules.
\end{enumerate}
\begin{rem}
Note that the tensor product $M\otimes_B-: \DD_{\sg}(B)\rightarrow \DD_{\sg}(A)$ is an equivalence of triangulated categories with the quasi-inverse $\Sigma^l (N\otimes_A-): \DD_{\sg}(A)\rightarrow \DD_{\sg}(B).$  Similarly, we have the following equivalence of triangulated categories
$$\Sigma^l (M\otimes_B-\otimes_BN): \DD_{\sg}(B\otimes B^{\op})\rightarrow \DD_{\sg}(A\otimes A^{\op}).$$
\end{rem}

Let us now prove the main result of this paper. 
\begin{thm}\label{thm}
Let $A$ and $B$ be two Noetherian  algebras over a field $k$ such that the enveloping algebras $A\otimes A^{\op}$ and $B\otimes B^{\op}$ are Noetherian. Suppose that $(_AM_B, _BN_A)$ defines a singular equivalence of Morita type with level $l\in\Z_{\geq 0}$. Then the  functor $\Sigma^l (M\otimes_B-\otimes_BN)$ induces an isomorphism of Gerstenhaber algebras between Tate-Hochschild cohomology rings $\HH_{\sg}^*(A, A)$ and $\HH_{\sg}^*(B, B)$.
\end{thm}
\begin{proof} First from the facts that the functor $\Sigma^l (M\otimes_B-\otimes_BN)$ induces an equivalence between $\DD_{\sg}(B\otimes B^{\op})$ and $\DD_{\sg}(A\otimes A^{\op})$ and that  the cup product $\cup$ in $\HH_{\sg}^*(A, A)$ coincides with the Yoneda product in $\DD_{\sg}(A\otimes A^{\op}),$ it follows that $\Sigma^l (M\otimes_B-\otimes_BN)$ yields an isomorphism of graded-commutative algebras between $\HH^*(B, B)$ and $\HH^*(A, A)$. It remains to prove that $\Sigma^l (M\otimes_B-\otimes_BN)$ induces an isomorphism of graded Lie algebras. For this,   $\Sigma^l (M\otimes_B-\otimes_BN)$ induces an isomorphism  between $\sgDPic_B$ and $\sgDPic_A$ and thus induces an isomorphism between $\Lie \sgDPic_B$ and $\Lie\sgDPic_A.$ In particular, this restricts to an isomorphism of graded Lie algebras $\Sigma^l (M\otimes_B-\otimes_BN): G_B\rightarrow G_A,$ where we denote $G_A:=\bigoplus_{m\in \Z}G_A(\epsilon_m).$ Consider the following commutative diagram
\begin{equation*}
\xymatrix@C=6pc{ G_B\ar[r]_-{\cong}^-{\Sigma^l (M\otimes_B-\otimes_BN)} & G_A\\\HH^{*+1}_{\sg}(B, B)\ar[u]^{\cong} \ar[r]_-{\cong}^-{\Sigma^l (M\otimes_B-\otimes_BN)} &   \HH^{*+1}_{\sg}(A, A)\ar[u]^-{\cong}.}
\end{equation*}
Since it follows from Corollary \ref{cor34} that the vertical morphisms are isomorphisms of graded Lie algebras,  the bottom  horizontal map induces an isomorphism of Gerstenhaber  algebras between $\HH_{\sg}^{*}(B, B)$ and $\HH^{*}_{\sg}(A, A)$. This proves the theorem. \end{proof}

\begin{cor}\label{cor-main}
Let $A$ and $B$ be two Noetherian $k$-algebras such that the enveloping algebras $A\otimes A^{\op}$ and $B\otimes B^{\op}$ are Noetherian. Assume that the derived categories $\DD(A)$ and $\DD(B)$ are equivalent as triangulated categories. Then there exists an isomorphism of Gerstenhaber algebras between $\HH_{\sg}^*(A, A)$ and $\HH_{\sg}^*(B, B)$. 
\end{cor}
\begin{proof} This comes from Theorem \ref{thm} and the fact that two derived equivalent algebras induce a singular equivalence of Morita type with some level $l\in\Z_{\geq 0}$ (cf. \cite{Wan14}).  \end{proof}

\bibliographystyle{plain}

\begin{thebibliography}{99}

\bibitem[BeLu]{BeLu}
Joseph Bernstein and Valery Lunts,
{\normalsize \it Equivariant sheaves and functors,}
Lecture Notes in Mathematics {\bf 1578},
Springer-Verlag, Berlin, 1994.


\bibitem[Buc]{Buc}
Ragnar-Olaf Buchweitz,
{\normalsize\it Maximal Cohen-Macaulay modules and Tate-cohomology over Gorenstein rings,} Universit\"at Hannover, 1986, available at http://hdl.handle.net/1807/16682.  



\bibitem[Del]{Deli}
Pierre Deligne,
{\normalsize\it SGA $4\frac{1}{2}$-Cohomologie \'etale,}
Lecture Notes in Mathematics {\bf 569}, Springer-Verlag, 1977.

%
%
%

\bibitem[Ger]{Ger1}
Murray Gerstenhaber,
{\normalsize\it The cohomology structure of an associative ring,}
Annals of Mathematics {\bf 78}(2) (1963),  267-288.

\bibitem[Kel98]{Kel2}
Bernhard Keller,
{\normalsize\it Invariance and localization for cyclic homology of DG algebras,}
Journal of Pure and Applied Algebra {\bf 123} (1998), 223-273.



\bibitem[Kel99]{Kel}
Bernhard Keller,
{\normalsize\it Hochschild cohomology and derived Picard groups,}
Journal of Pure and Applied Algebra {\bf 136} (1999), 1-56.

\bibitem[Kel18]{Kel18}
Bernhard Keller, 
{\normalsize\it Singular Hochschild cohomology via the singularity category,}
Comptes Rendus de l'Acad\'emie des Sciences - Series I - Mathematics {\bf 356}(11-12) (2018), 1106-1111. 



\bibitem[KeVo]{KeVo}
Bernhard Keller and Dieter Vossieck,
{\normalsize\it Sous les cat\'egories d\'eriv\'ees,} Comptes Rendus de l'Acad\'emie des Sciences - Series I - Mathematics {\bf 305}(6) (1987), 225-228.
\bibitem[KLZ]{KLZ}
Steffen Koenig, Yuming Liu, and Guodong Zhou, 
{\normalsize\it Transfer maps in Hochschild (co)homology and applications to stable and derived invariants and to the Auslander-Reiten conjecture,}
Transaction of the American Society {\bf 364} (1) (2012), 195-232. 



\bibitem[Lod]{Lod}
Jean-Louis Loday,
{\normalsize\it Cyclic homology,}
Grundlehren der Mathematischen Wissenschaften, 301, Springer, 1992.



\bibitem[Orl]{Orl}
Dmitri Orlov, 
{\normalsize\it Triangulated categories of singularities and D-branes in Landau-Ginzburg models,} Trudy Matematicheskogo Instituta Steklova {\bf 246} (2004), Algebr. Geom. Metody, Svyazi i Prilozh., 240-262. 

\bibitem[Ric]{Ric1}
Jeremy Rickard,
{\normalsize\it Morita theory for derived categories,}
Journal of London Mathematical Society {\bf 39}(2) (1989), 303-317.

\bibitem[Sta]{Sta}
Jim Stasheff,
{\normalsize\it The intrinsic bracket on the deformation complex
of an associative algebra,}
Journal of Pure and Applied Algebra {\bf 89} (1993), 231-235.

\bibitem[Wan15a]{Wan15a}
Zhengfang Wang, 
{\normalsize\it Singular Hochschild Cohomology and Gerstenhaber Algebra Structure,}
arXiv:1508.00190. 

\bibitem[Wan15b]{Wan14}
Zhengfang Wang,
{\normalsize\it Singular Equivalence of Morita Type with Level,}
Journal of Algebra {\bf 439} (2015), 245-269.



\bibitem[Wan18]{Wan15}
Zhengfang Wang,
{\normalsize\it Gerstenhaber algebra and Deligne's conjecture on Tate-Hochschild cohomology,}
arXiv:1801.07990. To appear in Transactions of the American Mathematical Society. 


\bibitem[Wei]{Wei}
Charles Weibel,
{\normalsize\it An introduction to Homological algebra,}
Cambridge University Press, Cambridge 1995.








\bibitem[Zim]{Zim}
Alexander Zimmermann,
{\normalsize\it Representation Theory: A Homological Algebra Point of View,}
Springer Verlag London, 2014.





\end{thebibliography}

\end{document}